\documentclass[fleqn,12pt,lettersize]{article}
\usepackage{amsthm}

\usepackage[latin1]{inputenc}
\usepackage{wasysym} % für \leadsto
\usepackage{fancyhdr}
\usepackage[all,dvips]{xy}
\usepackage[dvips,pdftitle={Locality for Classical Logic}, pdfauthor={Kai Brünnler}]{hyperref}

\def\contr  {{\css\mkern-0mu\tss}}%
\def\intr  {{\iss{\downarrow}}}%
\def\aintr {{\iss{\uparrow  }}}%
\let\swir=\sss%
\def\sysS {{\script S}}%
\def\merge {\mathbin{\bigdiamond}}%

\let\colder=<

\def\grammareq {\mathrel{\raise.4pt\hbox{::}{=}}}%

\def\permover{\mathchoice
             {\displaystyle
              \mathrel{\lower.493\fontdimen5\textfont3
                       \hbox{$\lhook$}%
                       \mkern-3mu{\rightarrow}}}%
             {\textstyle
              \mathrel{\lower.493\fontdimen5\textfont3
                       \hbox{$\lhook$}%
                       \mkern-3mu{\rightarrow}}}%
             {\scriptstyle
              \mathrel{\lower.54\fontdimen5\scriptfont3
                       \hbox{$\scriptstyle\lhook$}%
                       \mkern-3.2mu{\rightarrow}}}%
             {\scriptscriptstyle
              \mathrel{\lower.55\fontdimen5\scriptscriptfont3
                       \hbox{$\scriptscriptstyle\lhook$}%
                       \mkern-3.2mu{\rightarrow}}}%
             }%
\let\neg=\bar

\def\nsy #1#2#3#4{\nsysilent{#1}{#2}{#3}{#4}%
   \proofreadingcolor{$\csname #1\endcsname$}}%
\def\nsysilent #1#2#3#4{%
   \expandafter\let\csname #1\endcsname=#2%
   \sidx{#3@${#2}$, #4}}%

\def\boxtherule #1#2#3#4{\hbox{%
      \vtop{\kern0pt\hbox to0pt{\hss\structcolor{#1}\strut\enspace}}%
      \vtop{\kern0pt\structcolor{\boxit{\hbox to#2{\hfil\vbox to#3{%
                  \vfil\hbox to0pt{\hss\Black{$#4$}\hss}\vfil}\hfil}}}}}}%
\def\lowerdertobase #1{\lower
                       \ifsmallprint
                          6.89164pt
                       \else\ifverysmallprint
                          6.14815pt
                       \else
                          7.63518pt
                       \fi\fi\hbox{$#1$}}%
\def\InvisibleMark {\White{\vbox to0pt{\vss
   \hbox to0pt{\hss\vrule height1sp depth0pt width1sp}}}}%
\def\InvisibleMarkDown  #1{\kern-.#1pc\vbox to0pt{\kern.#1pc\InvisibleMark\vss}}%
\def\InvisibleMarkDDown #1{\kern-#1pc\vbox  to0pt{\kern#1pc\InvisibleMark \vss}}%
\def\InvisibleMarkUp    #1{\vbox to0pt{\vss\InvisibleMark\kern.#1pc}\kern-.#1pc}%
\def\InvisibleMarkUUp   #1{\vbox to0pt{\vss\InvisibleMark\kern#1pc}\kern-#1pc}%
\def\lv #1{\underline{#1\phantom{\hbox to0pt{,\hss}}}{}\mkern-1mu\lower1ex\hbox
                                                   {$\scriptscriptstyle\Vsss$}}%
\def\vl #1{\underline{#1\phantom{\hbox to0pt{,\hss}}}{}\lower1ex\hbox
                                                   {$\scriptscriptstyle\Lsss$}}%
\def\ls #1{\underline{#1\phantom{\hbox to0pt{,\hss}}}{}\lower1ex\hbox
                                                   {$\scriptscriptstyle\Ssss$}}%
\def\lg #1{\underline{#1\phantom{\hbox to0pt{,\hss}}}{}\lower1ex\hbox
                                                   {$\scriptscriptstyle\Gsss$}}%

\def\conldel {\{}%
\def\conrdel {\}}%
\def\lrgldel {\mathchoice{(}{(}{\langle}{\langle}}%
\def\lrgrdel {\mathchoice{)}{)}{\rangle}{\rangle}}%
\def\aprldel {\mathchoice
   {\mathopen {\setbox0=\hbox{$\displaystyle     \lrgldel$}\hbox to\wd0
                        {\hfil$\displaystyle     (       $\hfil}}}%
   {\mathopen {\setbox0=\hbox{$\textstyle        \lrgldel$}\hbox to\wd0
                        {\hfil$\textstyle        (        $\hfil}}}%
   {\mathopen {\setbox0=\hbox{$\scriptstyle      \lrgldel$}\hbox to\wd0
                        {\hfil$\scriptstyle      (        $\hfil}}}%
   {\mathopen {\setbox0=\hbox{$\scriptscriptstyle\lrgldel$}\hbox to\wd0
                        {\hfil$\scriptscriptstyle(        $\hfil}}}}%
\def\aprrdel {\mathchoice
   {\mathclose{\setbox0=\hbox{$\displaystyle     \lrgrdel$}\hbox to\wd0
                        {\hfil$\displaystyle     )       $\hfil}}}%
   {\mathclose{\setbox0=\hbox{$\textstyle        \lrgrdel$}\hbox to\wd0
                        {\hfil$\textstyle        )        $\hfil}}}%
   {\mathclose{\setbox0=\hbox{$\scriptstyle      \lrgrdel$}\hbox to\wd0
                        {\hfil$\scriptstyle      )        $\hfil}}}%
   {\mathclose{\setbox0=\hbox{$\scriptscriptstyle\lrgrdel$}\hbox to\wd0
                        {\hfil$\scriptscriptstyle)        $\hfil}}}}%
\def\seqldel {\mathchoice
   {\mathopen {\setbox0=\hbox{$\displaystyle     \lrgldel$}\hbox to\wd0
                        {\hfil$\displaystyle     \langle  $\hfil}}}%
   {\mathopen {\setbox0=\hbox{$\textstyle        \lrgldel$}\hbox to\wd0
                        {\hfil$\textstyle        \langle  $\hfil}}}%
   {\mathopen {\setbox0=\hbox{$\scriptstyle      \lrgldel$}\hbox to\wd0
                        {\hfil$\scriptstyle      \langle  $\hfil}}}%
   {\mathopen {\setbox0=\hbox{$\scriptscriptstyle\lrgldel$}\hbox to\wd0
                        {\hfil$\scriptscriptstyle\langle  $\hfil}}}}%
\def\seqrdel {\mathchoice
   {\mathclose{\setbox0=\hbox{$\displaystyle     \lrgrdel$}\hbox to\wd0
                        {\hfil$\displaystyle     \rangle  $\hfil}}}%
   {\mathclose{\setbox0=\hbox{$\textstyle        \lrgrdel$}\hbox to\wd0
                        {\hfil$\textstyle        \rangle  $\hfil}}}%
   {\mathclose{\setbox0=\hbox{$\scriptstyle      \lrgrdel$}\hbox to\wd0
                        {\hfil$\scriptstyle      \rangle  $\hfil}}}%
   {\mathclose{\setbox0=\hbox{$\scriptscriptstyle\lrgrdel$}\hbox to\wd0
                        {\hfil$\scriptscriptstyle\rangle  $\hfil}}}}%
\def\parldel {\mathchoice
   {\mathopen {\setbox0=\hbox{$\displaystyle     \lrgldel$}\hbox to\wd0
                        {\hfil$\displaystyle     [       $\hfil}}}%
   {\mathopen {\setbox0=\hbox{$\textstyle        \lrgldel$}\hbox to\wd0
                        {\hfil$\textstyle        [        $\hfil}}}%
   {\mathopen {\setbox0=\hbox{$\scriptstyle      \lrgldel$}\hbox to\wd0
                        {\hfil$\scriptstyle      [        $\hfil}}}%
   {\mathopen {\setbox0=\hbox{$\scriptscriptstyle\lrgldel$}\hbox to\wd0
                        {\hfil$\scriptscriptstyle[        $\hfil}}}}%
\def\parrdel {\mathchoice
   {\mathclose{\setbox0=\hbox{$\displaystyle     \lrgrdel$}\hbox to\wd0
                        {\hfil$\displaystyle     ]       $\hfil}}}%
   {\mathclose{\setbox0=\hbox{$\textstyle        \lrgrdel$}\hbox to\wd0
                        {\hfil$\textstyle        ]        $\hfil}}}%
   {\mathclose{\setbox0=\hbox{$\scriptstyle      \lrgrdel$}\hbox to\wd0
                        {\hfil$\scriptstyle      ]        $\hfil}}}%
   {\mathclose{\setbox0=\hbox{$\scriptscriptstyle\lrgrdel$}\hbox to\wd0
                        {\hfil$\scriptscriptstyle]        $\hfil}}}}%
\def\aprs #1{\aprldel #1\aprrdel}%
\def\cons #1{\conldel #1\conrdel}%
\def\pars #1{\parldel #1\parrdel}%
\def\sqn  #1{{\turnstile #1}}%

\def\quadcm {\rlap{\quad,}}%
\catcode`@=11
\def\upsmash{\relax % \relax, in case this comes first in \halign
  \ifmmode\def\next{\mathpalette\mathupsm@sh}\else\let\next\makeupsm@sh
  \fi\next}
\def\makeupsm@sh#1{\setbox\z@\hbox{#1}\finupsm@sh}
\def\mathupsm@sh#1#2{\setbox\z@\hbox{$\m@th#1{#2}$}\finupsm@sh}
\def\finupsm@sh{\ht\z@\z@ \box\z@}
\def\downsmash{\relax % \relax, in case this comes first in \halign
  \ifmmode\def\next{\mathpalette\mathdownsm@sh}\else\let\next\makedownsm@sh
  \fi\next}
\def\makedownsm@sh#1{\setbox\z@\hbox{#1}\findownsm@sh}
\def\mathdownsm@sh#1#2{\setbox\z@\hbox{$\m@th#1{#2}$}\findownsm@sh}
\def\findownsm@sh{\dp\z@\z@ \box\z@}
\catcode`@=12

\input rotate.tex

\def\hexnumber #1{\ifcase #10\or 1\or 2\or 3\or 4\or 5\or 6\or 7\or 8\or
   9\or A\or B\or C\or D\or E\or F\fi}%

\newfam\cmsmrfam
   \font\twelvesmr=cmsmr10 at 12pt
   \font\tensmr=cmsmr10
   \font\ninesmr=cmsmr10 at 9pt
   \font\eightsmr=cmsmr10 at 8pt
   \font\sevensmr=cmsmr7
   \font\sixsmr=cmsmr7 at 6pt
   \font\fivesmr=cmsmr5
   \textfont\cmsmrfam=\tensmr \scriptfont\cmsmrfam=\sevensmr
      \scriptscriptfont\cmsmrfam=\fivesmr

   \font\ninesmiu=cmsmiu10 at 9pt
   \font\eightsmiu=cmsmiu10 at 8pt
\newfam\eufmfam
   \font\teneufm=eufm10
   \font\nineeufm=eufm10 at 9pt
   \font\eighteufm=eufm10 at 8pt
   \font\seveneufm=eufm7
   \font\sixeufm=eufm7 at 6pt
   \font\fiveeufm=eufm5
   \textfont\eufmfam=\teneufm \scriptfont\eufmfam=\seveneufm
      \scriptscriptfont\eufmfam=\fiveeufm

\newfam\eurmfam
   \font\teneurm=eurm10
   \font\nineeurm=eurm10 at 9pt
   \font\eighteurm=eurm10 at 8pt
   \font\seveneurm=eurm7
   \font\sixeurm=eurm7 at 6pt
   \font\fiveeurm=eurm5
   \textfont\eurmfam=\teneurm \scriptfont\eurmfam=\seveneurm
      \scriptscriptfont\eurmfam=\fiveeurm

\newfam\eusmfam
   \font\teneusm=eusm10
   \font\nineeusm=eusm10 at 9pt
   \font\eighteusm=eusm10 at 8pt
   \font\seveneusm=eusm7
   \font\sixeusm=eusm7 at 6pt
   \font\fiveeusm=eusm5
   \textfont\eusmfam=\teneusm \scriptfont\eusmfam=\seveneusm
      \scriptscriptfont\eusmfam=\fiveeusm

\newfam\msamfam
   \font\tenmsam=msam10
   \font\ninemsam=msam10 at 9pt
   \font\eightmsam=msam10 at 8pt
   \font\sevenmsam=msam7
   \font\sixmsam=msam7 at 6pt
   \font\fivemsam=msam5
   \textfont\msamfam=\tenmsam \scriptfont\msamfam=\sevenmsam
      \scriptscriptfont\msamfam=\fivemsam

\newfam\msbmfam
   \font\tenmsbm=msbm10
   \font\ninemsbm=msbm10 at 9pt
   \font\eightmsbm=msbm10 at 8pt
   \font\sevenmsbm=msbm7
   \font\sixmsbm=msbm7 at 6pt
   \font\fivemsbm=msbm5
   \textfont\msbmfam=\tenmsbm \scriptfont\msbmfam=\sevenmsbm
      \scriptscriptfont\msbmfam=\fivemsbm

\newfam\scriptfam
      \font\tenfs=rsfs10 \skewchar\tenfs='177
      \font\ninefs=rsfs10 at 9pt \skewchar\ninefs='177
      \font\eightfs=rsfs10 at 8pt \skewchar\eightfs='177
      \font\sevenfs=rsfs7 \skewchar\sevenfs='177
      \font\sixfs=rsfs7 at 6pt \skewchar\sixfs='177
      \font\fivefs=rsfs5 \skewchar\fivefs='177
      \textfont\scriptfam=\tenfs \scriptfont\scriptfam=\sevenfs
         \scriptscriptfont\scriptfam=\fivefs
\newfam\stmaryrdfam
   \font\tenstmaryrd=stmary10
   \font\ninestmaryrd=stmary10 at 9pt
   \font\eightstmaryrd=stmary10 at 8pt
   \font\sevenstmaryrd=stmary7
   \font\sixstmaryrd=stmary7 at 6pt
   \font\fivestmaryrd=stmary5
   \textfont\stmaryrdfam=\tenstmaryrd \scriptfont\stmaryrdfam=\sevenstmaryrd
      \scriptscriptfont\stmaryrdfam=\fivestmaryrd

\def\eightpointlogic{%
   \textfont\eurmfam=\eighteurm \scriptfont\eurmfam=\sixeurm
      \scriptscriptfont\eurmfam=\fiveeurm
   \textfont\cmsmrfam=\eightsmr \scriptfont\cmsmrfam=\sixsmr
      \scriptscriptfont\cmsmrfam=\fivesmr
   \textfont\cmsmiufam=\eightsmiu \scriptfont\cmsmiufam=\sixsmiu
      \scriptscriptfont\cmsmiufam=\fivesmiu
   \textfont\eufmfam=\eighteufm \scriptfont\eufmfam=\sixeufm
      \scriptscriptfont\eufmfam=\fiveeufm
   \textfont\eusmfam=\eighteusm \scriptfont\eusmfam=\sixeusm
      \scriptscriptfont\eusmfam=\fiveeusm
   \textfont\msamfam=\eightmsam \scriptfont\msamfam=\sixmsam
      \scriptscriptfont\msamfam=\fivemsam
   \textfont\msbmfam=\eightmsbm \scriptfont\msbmfam=\sixmsbm
      \scriptscriptfont\msbmfam=\fivemsbm
   \ifscriptfamrsfs
      \textfont\scriptfam=\eightfs \scriptfont\scriptfam=\sixfs
         \scriptscriptfont\scriptfam=\fivefs
      \else
      \textfont\scriptfam=\eightmptwo \scriptfont\scriptfam=\sixmptwo
         \scriptscriptfont\scriptfam=\fivemptwo
      \fi
   \textfont\stmaryrdfam=\eightstmaryrd \scriptfont\stmaryrdfam=\sixstmaryrd
      \scriptscriptfont\stmaryrdfam=\fivestmaryrd}%
\let\oldeightpoint=\eightpoint
\def\eightpoint {\oldeightpoint\eightpointlogic}%

\def\ninepointlogic{%
   \textfont\eurmfam=\nineeurm \scriptfont\eurmfam=\sixeurm
      \scriptscriptfont\eurmfam=\fiveeurm
   \textfont\cmsmrfam=\ninesmr \scriptfont\cmsmrfam=\sixsmr
      \scriptscriptfont\cmsmrfam=\fivesmr
   \textfont\cmsmiufam=\ninesmiu \scriptfont\cmsmiufam=\sixsmiu
      \scriptscriptfont\cmsmiufam=\fivesmiu
   \textfont\eufmfam=\nineeufm \scriptfont\eufmfam=\sixeufm
      \scriptscriptfont\eufmfam=\fiveeufm
   \textfont\eusmfam=\nineeusm \scriptfont\eusmfam=\sixeusm
      \scriptscriptfont\eusmfam=\fiveeusm
   \textfont\msamfam=\ninemsam \scriptfont\msamfam=\sixmsam
      \scriptscriptfont\msamfam=\fivemsam
   \textfont\msbmfam=\ninemsbm \scriptfont\msbmfam=\sixmsbm
      \scriptscriptfont\msbmfam=\fivemsbm
   \ifscriptfamrsfs
      \textfont\scriptfam=\ninefs \scriptfont\scriptfam=\sixfs
         \scriptscriptfont\scriptfam=\fivefs
      \else
      \textfont\scriptfam=\ninemptwo \scriptfont\scriptfam=\sixmptwo
         \scriptscriptfont\scriptfam=\fivemptwo
      \fi
   \textfont\stmaryrdfam=\ninestmaryrd \scriptfont\stmaryrdfam=\sixstmaryrd
      \scriptscriptfont\stmaryrdfam=\fivestmaryrd}%
\let\oldninepoint=\ninepoint
\def\ninepoint {\oldninepoint\ninepointlogic}%

\def\twelvepointlogic{%
   \textfont\cmsmrfam=\twelvesmr \scriptfont\cmsmrfam=\eightsmr
      \scriptscriptfont\cmsmrfam=\sixsmr
      }%
\let\oldtwelvepoint=\twelvepoint
\def\twelvepoint {\oldtwelvepoint\twelvepointlogic}%

\def\script {\fam\scriptfam}%                        % Useless with zpmp2

\def\smaller #1{{\mathchoice{\scriptstyle{#1}}{\scriptstyle{#1}}%
      {\scriptscriptstyle{#1}}{\scriptscriptstyle{#1}}}}%

\catcode`\!=\active
   \edef!{\hexnumber\eurmfam}%
      \mathchardef\Deltabb="0!01
      \mathchardef\Pibb="0!05
      \mathchardef\Sigmabb="0!06
      \mathchardef\gammabb="0!0D
      \mathchardef\deltabb="0!0E
      \mathchardef\kappabb="0!14
      \mathchardef\pibb="0!19
      \mathchardef\psibb="0!20
      \mathchardef\rhobb="0!1A
      \mathchardef\sigmabb="0!1B
      \mathchardef\Abb="0!41
      \mathchardef\Bbb="0!42
      \mathchardef\Cbb="0!43
      \mathchardef\Dbb="0!44
      \mathchardef\Ebb="0!45
      \mathchardef\Fbb="0!46 
      \mathchardef\Gbb="0!47
      \mathchardef\Hbb="0!48
      \mathchardef\Ibb="0!49
      \mathchardef\Jbb="0!4A
      \mathchardef\Kbb="0!4B
      \mathchardef\Lbb="0!4C
      \mathchardef\Mbb="0!4D
      \mathchardef\Nbb="0!4E
      \mathchardef\Obb="0!4F
      \mathchardef\Pbb="0!50
      \mathchardef\Qbb="0!51
      \mathchardef\Rbb="0!52
      \mathchardef\Sbb="0!53
      \mathchardef\Tbb="0!54
      \mathchardef\Ubb="0!55
      \mathchardef\Vbb="0!56
      \mathchardef\Wbb="0!57
      \mathchardef\Xbb="0!58
      \mathchardef\Ybb="0!59
      \mathchardef\Zbb="0!5A
      \mathchardef\abb="0!61
      \mathchardef\bbb="0!62
      \mathchardef\cbb="0!63
      \mathchardef\dbb="0!64
      \mathchardef\ebb="0!65
      \mathchardef\fbb="0!66
      \mathchardef\gbb="0!67
      \mathchardef\hbb="0!68
      \mathchardef\ibb="0!69
      \mathchardef\jbb="0!6A
      \mathchardef\kbb="0!6B
      \mathchardef\lbb="0!6C
      \mathchardef\mbb="0!6D
      \mathchardef\nbb="0!6E
      \mathchardef\obb="0!6F
      \mathchardef\pbb="0!70
      \mathchardef\qbb="0!71
      \mathchardef\rbb="0!72
      \mathchardef\sbb="0!73
      \mathchardef\tbb="0!74
      \mathchardef\ubb="0!75
      \mathchardef\vbb="0!76
      \mathchardef\wbb="0!77
      \mathchardef\xbb="0!78
      \mathchardef\ybb="0!79
      \mathchardef\zbb="0!7A
   \edef!{\hexnumber\cmsmrfam}%
      \mathchardef\Ass="0!41
      \mathchardef\Bss="0!42
      \mathchardef\Css="0!43
      \mathchardef\Dss="0!44
      \mathchardef\Ess="0!45
      \mathchardef\Fss="0!46 
      \mathchardef\Gss="0!47
      \mathchardef\Hss="0!48
      \mathchardef\Iss="0!49
      \mathchardef\Jss="0!4A
      \mathchardef\Kss="0!4B
      \mathchardef\Lss="0!4C
      \mathchardef\Mss="0!4D
      \mathchardef\Nss="0!4E
      \mathchardef\Oss="0!4F
      \mathchardef\Pss="0!50
      \mathchardef\Qss="0!51
      \mathchardef\Rss="0!52
      \mathchardef\Sss="0!53
      \mathchardef\Tss="0!54
      \mathchardef\Uss="0!55
      \mathchardef\Vss="0!56
      \mathchardef\Wss="0!57
      \mathchardef\Xss="0!58
      \mathchardef\Yss="0!59
      \mathchardef\Zss="0!5A
      \mathchardef\ass="0!61
      \mathchardef\bss="0!62
      \mathchardef\css="0!63
      \mathchardef\dss="0!64
      \mathchardef\ess="0!65
      \mathchardef\fss="0!66
      \mathchardef\gss="0!67
      \mathchardef\hhss="0!68    % \hss is already defined!!!
      \mathchardef\iss="0!69
      \mathchardef\jss="0!6A
      \mathchardef\kss="0!6B
      \mathchardef\lss="0!6C
      \mathchardef\mss="0!6D
      \mathchardef\nss="0!6E
      \mathchardef\oss="0!6F
      \mathchardef\pss="0!70
      \mathchardef\qss="0!71
      \mathchardef\rss="0!72
      \mathchardef\sss="0!73
      \mathchardef\tss="0!74
      \mathchardef\uss="0!75
      \mathchardef\vvss="0!76    % \vss is already defined!!!
      \mathchardef\wss="0!77
      \mathchardef\xss="0!78
      \mathchardef\yss="0!79
      \mathchardef\zss="0!7A
      \mathchardef\bang="0!21
      \mathchardef\loc="0!21
      \mathchardef\lone="0!31
      \mathchardef\ltwo="0!32   % kai: test
      \mathchardef\lthree="0!33
      \mathchardef\lwn="0!3F
      \mathchardef\lzero="0!30
      \mathchardef\Gammass="0!00
      \mathchardef\Deltass="0!01
      \mathchardef\deltass="0!0E
      \mathchardef\gitss="0!67
   \edef!{\hexnumber\msamfam}%
      \mathchardef\bigdiamond="0!06
      \mathchardef\triangleup="0!4D
      \mathchardef\triangledown="0!4F
   \edef!{\hexnumber\msbmfam}%
      \mathchardef\emptyset="0!3F
   \edef!{\hexnumber\scriptfam}%
         \mathchardef\Asc="0!41
         \mathchardef\Bsc="0!42
         \mathchardef\Csc="0!43
         \mathchardef\Dsc="0!44
         \mathchardef\Esc="0!45
         \mathchardef\Fsc="0!46 
         \mathchardef\Gsc="0!47
         \mathchardef\Hsc="0!48
         \mathchardef\Isc="0!49
         \mathchardef\Jsc="0!4A
         \mathchardef\Ksc="0!4B
         \mathchardef\Lsc="0!4C
         \mathchardef\Msc="0!4D
         \mathchardef\Nsc="0!4E
         \mathchardef\Osc="0!4F
         \mathchardef\Psc="0!50
         \mathchardef\Qsc="0!51
         \mathchardef\Rsc="0!52
         \mathchardef\Ssc="0!53
         \mathchardef\Tsc="0!54
         \mathchardef\Usc="0!55
         \mathchardef\Vsc="0!56
         \mathchardef\Wsc="0!57
         \mathchardef\Xsc="0!58
         \mathchardef\Ysc="0!59
         \mathchardef\Zsc="0!5A
   \edef!{\hexnumber\cmsmrfam}%
\catcode`\!=12
\mathchardef\contr="013E
\let\false=\Fss

\newbox\Nablassbox
\def\Nablass {{\mathchoice
      {\setbox\Nablassbox=\hbox{$\Deltass$}%
         \setbox\Nablassbox=\hbox{\rotu\Nablassbox}%
         \box\Nablassbox}%
      {\setbox\Nablassbox=\hbox{$\Deltass$}%
         \setbox\Nablassbox=\hbox{\rotu\Nablassbox}%
         \box\Nablassbox}%
      {\setbox\Nablassbox=\hbox{$\scriptstyle\Deltass$}%
         \setbox\Nablassbox=\hbox{\rotu\Nablassbox}%
         \box\Nablassbox}%
      {\setbox\Nablassbox=\hbox{$\scriptscriptstyle\Deltass$}%
         \setbox\Nablassbox=\hbox{\rotu\Nablassbox}%
         \box\Nablassbox}}}%

\let\true=\Tss
\mathchardef\weak="013C
\def\Gsss {{\smaller\Gss}}%
\def\Lsss {{\smaller\Lss}}%
\def\Ssss {{\smaller\Sss}}%
\def\Vsss {{\smaller\Vss}}%
\mathchardef\impl="221B
\catcode`\!=\active
   \edef!{\hexnumber\stmaryrdfam}%
   \mathchardef\squarebox="2!1F
   \mathchardef\lpar="2!4F
   \mathchardef\lplus="2!16
   \mathchardef\lmix="2!22
   \mathchardef\lprec="2!34
   \mathchardef\ltens="2!0F
   \mathchardef\lwith="2!4E
   \mathchardef\merge="2!05
   \edef!{\hexnumber\msamfam}%
   \mathchardef\limpalt="2!28
   \catcode`\!=12

\catcode`\!=\active
   \edef!{\hexnumber\msamfam}%
   \mathchardef\ge="3!3E
   \mathchardef\le="3!36
   \mathchardef\gex="3!3C
   \mathchardef\lex="3!34
   \edef!{\hexnumber\stmaryrdfam}%
   \mathchardef\bigcross="3!22
   \mathchardef\arrowequiv="3!2D
   \mathchardef\inplus="3!41
   \mathchardef\msin="3!41
   \mathchardef\subsetplus="3!44
   \mathchardef\supsetplus="3!45
   \mathchardef\subseteqplus="3!46
   \mathchardef\supseteqplus="3!47
   \catcode`\!=12

\let\turnstile=\vdash

\catcode`\!=\active
   \edef!{\hexnumber\stmaryrdfam}%
   \mathchardef\lbbrack="4!4A
   \mathchardef\rbbrack="5!4B
   \mathchardef\lbpar="4!4C
   \mathchardef\rbpar="5!4D
      \mathchardef\lstrange="4!2A
      \mathchardef\rstrange="5!2B
      \mathchardef\lstrange="4!48
      \mathchardef\rstrange="5!49
   \catcode`\!=12

\def\xyvdots {\raise6pt\hbox{$\vdots$}}%

\newdimen\dercldim                                % dcl
\newdimen\derccdim                                % dcc
\newdimen\dercrdim                                % dcr
\newdimen\derldim                                 % dl
\newdimen\dercdim                                 % dc
\newdimen\derrdim                                 % dr
\newdimen\derdim                                  % d
\newdimen\derdldim                                % ddl
\newdimen\derdrdim                                % ddr
\newbox\derboxone                                 % b1
\newbox\derboxtwo                                 % b2
\newbox\derboxthree                               % b3
\newbox\derboxfour                                % b4
\newdimen\derquad\derquad=\fontdimen6\textfont2
\newdimen\deropen\deropen=\fontdimen5\textfont2\divide\deropen by3
\def\leaf #1{\global\setbox\derboxone=\hbox{\strut$#1$}%
   \global\derldim=0pt                            % dl=0
   \global\dercdim=\wd\derboxone                  % dc=wd(b1)
   \global\derrdim=0pt                            % dr=0
   }%
\def\rootaux #1#2#3{\setbox\derboxtwo=\hbox{\unhbox\derboxone}%
   \setbox\derboxthree=\hbox 
      {$\smash{\lower\fontdimen22\textfont2\hbox{$#1$}}$}%
   \setbox\derboxfour=\hbox 
      {$\smash{\lower\fontdimen22\textfont2\hbox{$#2$}}$}%
   \leaf{#3}%                                     % dl=0, dc=wd(b1), dr=0
   \derdim=\dercdim\advance\derdim by-\derccdim\divide\derdim by2 
   \global\derldim=\dercldim\global\advance\derldim by-\derdim
   \global\derrdim=\dercrdim\global\advance\derrdim by-\derdim
   \deropen=\fontdimen5\textfont2\divide\deropen by3
   \setbox\derboxone=\hbox{\vbox{\offinterlineskip
         \hbox{\ifdim\derldim<0pt\kern-\derldim\fi
               \box\derboxtwo
               \ifdim\derrdim<0pt\kern-\derrdim\fi}%
         \kern\deropen
         \hbox{\ifdim\dercldim>\derldim
                  \ifdim\derldim>0pt\kern\derldim\fi
                  \else\kern\dercldim\fi
               \hbox to0pt{\hss\copy\derboxthree}%
               \vbox{\ifdim\derccdim>\dercdim\hsize=\derccdim
                                        \else\hsize=\dercdim \fi
                    \hrule height.2pt depth.2pt width\hsize}%
               \hbox to0pt{\copy\derboxfour\hss}%
               \ifdim\dercrdim>\derrdim
                  \ifdim\derrdim>0pt\kern\derrdim\fi
                  \else\kern\dercrdim\fi}%
         \kern\deropen
         \hbox{\ifdim\derldim>0pt\kern\derldim\fi
               \box\derboxone
               \ifdim\derrdim>0pt\kern\derrdim\fi}}}%
   \ifdim\derldim<0pt\global\derldim=0pt\fi       % dl=max(dl,0)
   \ifdim\derrdim<0pt\global\derrdim=0pt\fi       % dr=max(dr,0)
   \derdldim=\wd\derboxthree\advance\derdldim by-\dercldim
   \derdrdim=\wd\derboxfour \advance\derdrdim by-\dercrdim
   \ifdim\derdim<0pt
      \ifdim\derdldim<0pt
         \derdldim=0pt                            % d<0, ddl<0 -> ddl=0
      \fi
      \ifdim\derdrdim<0pt
         \derdrdim=0pt                            % d<0, ddr<0 -> ddr=0
      \fi
   \else
      \ifdim\derldim>0pt
         \ifdim\derdldim>-\derdim
            \advance\derdldim by\derdim           % d>=0, dl>0, ddl+d>0 -> 
         \else                                            %                    ddl=ddl+d
            \derdldim=0pt                         % d>=0, dl>0, ddl+d<=0 ->
         \fi                                      %                        ddl=0
      \else
         \advance\derdldim by\dercldim            % d>=0, dl=0 -> ddl=ddl+dcl
      \fi
      \ifdim\derrdim>0pt
         \ifdim\derdrdim>-\derdim
            \advance\derdrdim by\derdim           % d>=0, dr>0, ddr+d>0 -> 
         \else                                            %                    ddr=ddr+d
            \derdrdim=0pt                         % d>=0, dr>0, ddr+d<=0 ->
         \fi                                      %                        ddr=0
      \else
         \advance\derdrdim by\dercrdim            % d>=0, dr=0 -> ddr=ddr+dcr
      \fi
   \fi
   \global\setbox\derboxone=\hbox
      {\kern\derdldim\unhbox\derboxone\kern\derdrdim}%
   \global\advance\derldim by\derdldim            % dl=dl+ddl
   \global\advance\derrdim by\derdrdim            % dr=dr+ddr
   }%
\def\rootr #1#2#3#4{{#4}%
   \dercldim=\derldim
   \derccdim=\dercdim
   \dercrdim=\derrdim
   \rootaux{#1}{#2}{#3}}%
\def\rrootr #1#2#3#4#5{\derquad=\fontdimen6\textfont2
   {#4}%
           \dercldim  =\derldim
   \setbox\derboxtwo=\hbox{\unhbox\derboxone\kern\derquad}%
           \derccdim  =\dercdim
   \advance\derccdim by\derrdim
   \advance\derccdim by\derquad
   {#5}%
   \setbox\derboxone=\hbox{\unhbox\derboxtwo\unhbox\derboxone}%
   \advance\derccdim by\derldim
   \advance\derccdim by\dercdim
           \dercrdim  =\derrdim
   \rootaux{#1}{#2}{#3}}%
\def\rrrootr #1#2#3#4#5#6{\derquad=\fontdimen6\textfont2
   {#4}%
           \dercldim  =\derldim
   \setbox\derboxtwo=\hbox{\unhbox\derboxone\kern\derquad}%
           \derccdim  =\dercdim
   \advance\derccdim by\derrdim
   \advance\derccdim by\derquad
   {#5}%
   \setbox\derboxtwo=\hbox{\unhbox\derboxtwo\unhbox\derboxone\kern\derquad}%
   \advance\derccdim by\derldim
   \advance\derccdim by\dercdim
   \advance\derccdim by\derrdim
   \advance\derccdim by\derquad
   {#6}%
   \setbox\derboxone=\hbox{\unhbox\derboxtwo\unhbox\derboxone}%
   \advance\derccdim by\derldim
   \advance\derccdim by\dercdim
           \dercrdim  =\derrdim
   \rootaux{#1}{#2}{#3}}%
\def\root       #1#2#3{\rootr  {#1\;}{}{#2}{#3}}%
\def\rootnote #1#2#3#4{\rootr {#1\;}{\;#2}{#3}{#4}}%
\def\rroot    #1#2#3#4{\rrootr {#1\;}{}{#2}{#3}{#4}}%
\def\deraux {\derldim=0pt\dercdim=0pt\derrdim=0pt}%
\def\der       #1#2#3{\deraux\root  {#1}{#2}{#3}        \box\derboxone}%
\def\dder    #1#2#3#4{\deraux\rroot {#1}{#2}{#3}{#4}    \box\derboxone}%
\def\dernote       #1#2#3#4{\deraux\rootr  {#1\;}{\;#2}{#3}{#4}\box\derboxone}%
\def\ddernote    #1#2#3#4#5{\deraux\rrootr {#1\;}{\;#2}{#3}{#4}{#5}\box
                                                                   \derboxone}%
\def\inf       #1#2#3{\der  {#1}{#2}{\leaf{#3}}}%
\def\iinf    #1#2#3#4{\dder {#1}{#2}{\leaf{#3}}{\leaf{#4}}}%
\def\infnote       #1#2#3#4{\dernote  {#1}{#4}{#2}{\leaf{#3}}}%
\def\iinfnote    #1#2#3#4#5{\ddernote {#1}{#5}{#2}{\leaf{#3}}{\leaf{#4}}}%
\newbox\derskelboxone
\newbox\derskelboxtwo
\newbox\derskelboxthree
\newbox\derskelboxfour
\newdimen\derskeldimenone
\newdimen\derskeldimentwo
\newdimen\derskeldimenthree
\newdimen\derskeldimenfour
\newdimen\derskeldimenfive
\newdimen\derskeldimensix
\newdimen\derskeldimenseven
\newdimen\derskeldimeneight
\def\derskel #1#2#3#4{%
   \setbox\derskelboxone=\hbox{$#1$\strut}%
   \derskeldimenone=\ht\derskelboxone
   \advance\derskeldimenone by\dp\derskelboxone
   \derskeldimentwo=\wd\derskelboxone
   \divide\derskeldimentwo by2
   \setbox\derskelboxone=\hbox to0pt{%
      \hss\raise\dp\derskelboxone\box\derskelboxone\hss}%
   \ht\derskelboxone=0pt
   \dp\derskelboxone=0pt
   \setbox\derskelboxtwo=\hbox{$#3$\strut}%
   \derskeldimenthree=\ht\derskelboxtwo
   \advance\derskeldimenthree by\dp\derskelboxtwo
   \derskeldimenfour=\wd\derskelboxtwo
   \divide\derskeldimenfour by2
   \setbox\derskelboxtwo=\hbox to0pt{%
      \hss\raise\dp\derskelboxtwo\box\derskelboxtwo\hss}%
   \ht\derskelboxtwo=0pt
   \dp\derskelboxtwo=0pt
   \ifdim\derskeldimenone>\derskeldimenthree
      \else\derskeldimenone=\derskeldimenthree\fi
   \setbox\derskelboxthree=\hbox{$#4$\strut}%
   \derskeldimenfive=\ht\derskelboxthree
   \advance\derskeldimenfive by\dp\derskelboxthree
   \derskeldimensix=\wd\derskelboxthree
   \divide\derskeldimensix by2
   \setbox\derskelboxthree=\hbox to0pt{%
      \hss\lower\ht\derskelboxthree\box\derskelboxthree\hss}%
   \ht\derskelboxthree=0pt
   \dp\derskelboxthree=0pt
   \setbox\derskelboxfour=\hbox{$#2$\strut}%
   \derskeldimenseven=\ht\derskelboxfour
   \advance\derskeldimenseven by\dp\derskelboxfour
   \derskeldimeneight=\wd\derskelboxfour
   \divide\derskeldimeneight by2
   \setbox\derskelboxfour=\hbox to0pt{%
      \hss\raise\dp\derskelboxfour\box\derskelboxfour\hss}%
   \ht\derskelboxfour=0pt
   \dp\derskelboxfour=0pt
   \ifdim\derskeldimenone>\derskeldimenseven
      \else\derskeldimenone=\derskeldimenseven\fi
   \derskeldimenthree=\derskeldimentwo
   \advance\derskeldimenthree by2\derskeldimeneight
   \advance\derskeldimenthree by\derskeldimenfour
   \advance\derskeldimenthree by2em
   \divide\derskeldimenthree by2
   \advance\derskeldimensix by-\derskeldimenthree
   \derskeldimenseven=\derskeldimensix
   \advance\derskeldimensix by-\derskeldimentwo
   \advance\derskeldimenseven by-\derskeldimenfour
   \ifdim\derskeldimensix>0pt
      \else\derskeldimensix=0pt\fi
   \ifdim\derskeldimenseven>0pt
      \else\derskeldimenseven=0pt\fi
   \vbox{\kern\derskeldimenone\hbox{\kern\derskeldimensix
         \kern\derskeldimentwo
         \xy
         <-\derskeldimenthree,\derskeldimenthree>="here"
            *{\box\derskelboxone}**\dir{-};
         "here"+<\derskeldimentwo,0pt>="here"**\dir{-};
         "here"+<1em,0pt>="here"**\dir{-};
         "here"+<\derskeldimeneight,0pt>="here"
            *{\box\derskelboxfour}**\dir{-};
         "here"+<\derskeldimeneight,0pt>="here"**\dir{-};
         "here"+<1em,0pt>="here"**\dir{-};
         "here"+<\derskeldimenfour,0pt>*{\box\derskelboxtwo}**\dir{-};
         0*{\box\derskelboxthree}**\dir{-};
         <-\derskeldimenthree,\derskeldimenthree>**\dir{-}
         \endxy
         \kern\derskeldimenfour\kern\derskeldimenseven}%
      \kern\derskeldimenfive}}%

\newbox\DerivOneBox
\newbox\DerivTwoBox
\newbox\DerivThreeBox
\newbox\DerivFourBox
\newdimen\DerivOneDimen
\newdimen\DerivTwoDimen
\newdimen\DerivThreeDimen
\newdimen\DerivFourDimen
\def\Derivation #1#2#3#4#5{\DerivationFactors{#1}{#2}{#3}{#4}{#5}11}%
\def\Derivationleaf #1#2#3#4#5{\global\setbox\derboxone=\hbox{\strut
                                    $\DerivationFactors{#1}{#2}{#3}{#4}{#5}11$}}%
\def\DerivationFactors #1#2#3#4#5#6#7{%
   \setbox\DerivOneBox=\hbox{$#1\strut$}%
      \DerivOneDimen=\wd\DerivOneBox\divide\DerivOneDimen by2
   \setbox\DerivThreeBox=\hbox{$#3\strut$}%
      \DerivThreeDimen=\wd\DerivThreeBox\divide\DerivThreeDimen by2
   \setbox\DerivTwoBox=\hbox{\box\DerivOneBox\hbox{$#2$}\box\DerivThreeBox}%
      \DerivTwoDimen=\wd\DerivTwoBox
   \setbox\DerivFourBox=\hbox{$#4\strut$}%
      \DerivFourDimen=\wd\DerivFourBox
   \ifdim\DerivFourDimen>\DerivTwoDimen
      \global\dercdim=\DerivFourDimen                % dc=wd(b4) see logicmac.tex
      \global\derldim=0pt                            % dl
      \global\derrdim=0pt                            % dr
      \advance\DerivFourDimen by-\DerivTwoDimen
      \divide \DerivFourDimen by2
      \advance\DerivTwoDimen  by-\DerivOneDimen
      \advance\DerivTwoDimen  by-\DerivThreeDimen
      \divide \DerivTwoDimen  by 2
   \else
      \global\dercdim=\DerivFourDimen                % dc=wd(b4) see logicmac.tex
      \DerivFourDimen=0pt
      \advance\DerivTwoDimen  by-\DerivOneDimen
      \advance\DerivTwoDimen  by-\DerivThreeDimen
      \global\derldim=\DerivTwoDimen
         \global\advance\derldim by-\dercdim
         \global\divide\derldim by2
         \global\advance\derldim by\DerivOneDimen    % dl
      \global\derrdim=\DerivTwoDimen
         \global\advance\derrdim by-\dercdim
         \global\divide\derrdim by2
         \global\advance\derrdim by\DerivThreeDimen  % dr
      \divide \DerivTwoDimen  by 2
   \fi
   \vbox{\offinterlineskip\hbox{\kern\DerivFourDimen\box\DerivTwoBox}%
         \hbox{\kern\DerivFourDimen\kern\DerivOneDimen
               \kern\DerivTwoDimen\kern-#6\DerivTwoDimen\hbox{$\xy
               0;<#6\DerivTwoDimen,0pt>:<0pt,#7\DerivTwoDimen>::
               (0,1);(2,1)**\crv{(1.25,1.1875)&(0.75,0.8125)};
               (1,0)**@{-};(0,1)**@{-};
               (1,0.625)*{\scriptstyle #5}
               \endxy$}}%
         \hbox{\kern\DerivFourDimen\kern\DerivOneDimen\kern\DerivTwoDimen
               \hbox to0pt{\hss\box\DerivFourBox\hss}%
               \kern\DerivFourDimen\kern\DerivOneDimen\kern\DerivTwoDimen}}}%
\newcommand{\ie}{i.e.\ }
\newcommand{\eg}{e.g.\ }
\newcommand{\cf}{cf.\ }

\newcommand{\smalltitle}[1]{{\bf #1}}
\newcommand{\vciinf}[4]{\vcenter{\iinf{#1}{#2}{#3}{#4}}}
  
\newcommand{\vcinf}[3]{\vcenter{\inf{#1}{#2}{#3}}}

\newcommand{\LK}{\sf{LK}}
\newcommand{\KS}{\sf{KS}}
\newcommand{\SKS}{\sf{SKS}}
\newcommand{\SKSpred}{\sf{SKSq}}
\newcommand{\KSpred}{\sf{KSq}}
\newcommand{\SKSg}{\sf{SKSg}}
\newcommand{\KSg}{\sf{KSg}}
\newcommand{\SKSgpred}{\sf{SKSgq}}
\newcommand{\KSgpred}{\sf{KSgq}}

\newcommand{\GSonep}{\sf{GS1p}}
\newcommand{\GSone}{\sf{GS1}}

\newcommand{\Gthree}{\sf{G3}}

\def\sruleAx {{\Ass\xss}}%
\def\sruleRW {{\Rss\Wss}}%
\def\sruleRC {{\Rss\Css}}%
\def\sruleCut {{\Css\uss\tss}}%
\def\sruleRv {{\Rss\vee}}%
\def\sruleRand {{\Rss\wedge}}%
\def\sruleRE {{\Rss\exists}}%
\def\sruleRA {{\Rss\forall}}%

\let\swir=\sss%
\def\rulem {{\mss}}%
\def\rules {{\sss}}%
\def\ruleacdown {{\ass\mkern-1.5mu\css{\downarrow  }}}%
\def\ruleacup {{\ass\mkern-1.5mu\css{\uparrow  }}}%
\def\ruleawdown {{\ass\mkern-1.5mu\wss{\downarrow  }}}%
\def\ruleawup {{\ass\mkern-1.5mu\wss{\uparrow  }}}%
\def\ruleaisdown {{\ass\mkern-1.5mu\iss_\sss{\downarrow  }}}%
\def\ruleaisup {{\ass\mkern-1.5mu\iss_\sss{\uparrow  }}}%
\def\ruleaidown {{\ass\mkern-1.5mu\iss{\downarrow  }}}%
\def\ruleaiup {{\ass\mkern-1.5mu\iss{\uparrow  }}}%
\def\rulessdown {{\sss\mkern-1mu\sss{\downarrow  }}}%
\def\rulessup {{\sss\mkern-1mu\sss{\uparrow  }}}%

\def\rulecdown {{\css{\downarrow  }}}%
\def\ruleddown {{\dss{\downarrow  }}}%
\def\rulecup {{\css{\uparrow  }}}%
\def\rulewdown {{\wss{\downarrow  }}}%
\def\rulewup {{\wss{\uparrow  }}}%
\def\ruleidown {{\iss{\downarrow  }}}%
\def\ruleiup {{\iss{\uparrow  }}}%

\def\true {{\sf t}}%
\def\false {{\sf f}}%

\def\rulendown {{\nss{\downarrow  }}}%
\def\rulenup {{\nss{\uparrow  }}}%
\def\ruleudown {{\uss{\downarrow  }}}%
\def\ruleuup {{\uss{\uparrow  }}}%

\def\rulemonedown {\sf{l_1}{\downarrow  }}%
\def\rulemoneup {\sf{l_1}{\uparrow  }}%
\def\rulemtwodown {\sf{l_2}{\downarrow  }}%
\def\rulemtwoup {\sf{l_2}{\uparrow  }}%

\newcommand{\subst}{\pars}
\newcommand{\tinyspace}{\:\!} % 1/18 von einem quad
\newcommand{\theoremnl}{\hspace{0mm}\\}

\begin{document}

\begin{titlepage}
%\thispagestyle{fancy}
%\lhead{\footnotesize {\textsc{Draft -- \today}} }
%%\lhead{}
%\rhead{}
%\lfoot{}
%\cfoot{}
%\rfoot{}
%\renewcommand{\headrulewidth}{0pt}
%\renewcommand{\footrulewidth}{0pt}

\textbf{\Large Locality for Classical Logic}  \\

  Kai Brünnler\\
  {\footnotesize Technische Universität Dresden\\
  Fakultät Informatik - 
  01062 Dresden -
  Germany\\
  \texttt{kai.bruennler@inf.tu-dresden.de}}\\

\vfill  
 
\begin{center}
\begin{minipage}[h]{12cm}
  \textbf{\small Abstract} \hspace{1ex}\footnotesize In this paper we
  will see deductive systems for classical propositional and predicate
  logic in the calculus of structures.  Like sequent systems, they
  have a cut rule which is admissible. In addition, they enjoy a
  top-down symmetry and some normal forms for derivations that are not
  available in the sequent calculus. Identity axiom, cut, weakening
  and also contraction can be reduced to atomic form.  This leads to
  rules that are \emph{local}: they do not require the inspection of
  expressions of unbounded size.
  \end{minipage}
\end{center}

\vfill
 \begin{center}
 \begin{minipage}[h]{12cm}
   \renewcommand{\contentsname}{\small Table of Contents}
   \footnotesize
\tableofcontents
\end{minipage}
\end{center}

\end{titlepage}

\theoremstyle{definition}
\newtheorem{Theorem}{Theorem}[section]
\newtheorem{Corollary}[Theorem]{Corollary}
\newtheorem{Proposition}[Theorem]{Proposition}
\newtheorem{Definition}[Theorem]{Definition}
\newtheorem{Lemma}[Theorem]{Lemma}
\newtheorem{Conjecture}[Theorem]{Conjecture}
\newtheorem{Remark}[Theorem]{Remark}
\newtheorem{Example}[Theorem]{Example}
\newtheorem{Problem}[Theorem]{Problem}
\newtheorem{Notation}[Theorem]{Notation}

\section{Introduction}

Inference rules that copy an unbounded quantity of information are
problematic from the points of view of complexity and implementation.
In the sequent calculus, an example is given by the contraction rule
in Gentzen's $\LK$ \cite{GenILD35}:
$$
\infnote{}{\Gamma \vdash \Phi, A}{\Gamma \vdash \Phi, A, A}{\quad .}
$$
Here, going from bottom to top in constructing a proof, a formula
$A$ of unbounded size is duplicated. Whatever mechanism performs this
duplication, it has to inspect all of $A$, so it has to have a
\emph{global} view on $A$. If, for example, we had to implement
contraction on a distributed system, where each processor has a
limited amount of local memory, the formula $A$ could be spread over a
number of processors. In that case, no single processor has a global
view on $A$, and we should put in place complex mechanisms to cope
with the situation. 

Let us call \emph{local} those inference rules that do not require
such a global view on formulae of unbounded size, and \emph{non-local}
those rules that do. Further examples of non-local rules are the
promotion rule in the sequent calculus for linear logic (left,
\cite{GirTCS87}) and context-sharing (or additive) rules found in
various sequent systems (right, \cite{TroSch96}):
$$
\vcenter{\infnote{} {\sqn{\loc A,\lwn B_1,\dots,\lwn B_n}} {\sqn{A ,\lwn
    B_1,\dots,\lwn B_n}}{}}\qquad \mbox{and}\qquad  
\vcenter{\iinfnote{}{\Gamma \vdash \Phi, A
  \wedge B}{\Gamma \vdash \Phi, A}{\Gamma \vdash \Phi, B}{\quad
  .}}
$$
To apply the promotion rule, one has to check whether all formulae
in the context are prefixed with a question mark modality: the number of
formulae to check is unbounded. To apply the context-sharing
$\sruleRand$ rule, a context of unbounded size has to be copied.

While there are methods to solve these problems in an implementation,
an interesting question is whether it is possible to approach them
proof-theoretically, \ie by avoiding non-local rules. The present work
gives an affirmative answer by presenting systems for both classical
propositional and first-order predicate logic in which context-sharing
rules as well as contraction are replaced by local rules. For
propositional logic it is even possible to obtain a system which
contains local rules only (which has already been presented in
\cite{BruTiu01}).

Locality is achieved by reducing the problematic rules to their atomic
forms. This is not entirely new: in most sequent systems for classical
logic the identity axiom is reduced to its atomic form, \ie
$$
\inf{}{A \vdash  A}{}
\qquad \mbox{is equivalently replaced by} \qquad
\inf{}{ a \vdash a}{}\quad,
$$
where $a$ is an atom. Contraction, however, cannot be replaced by
its atomic form in known sequent systems \cite{BruRC02}. In fact, I
believe that such a system cannot be presented in the sequent
calculus. To obtain local inference rules, I employ the \emph{calculus
  of structures} \cite{Gug02,GugStr01}.  This formalism differs from
the sequent calculus in two main aspects:
\begin{enumerate}
\item \emph{Deepness}: inference rules apply anywhere deep inside a
  formula, not only at the main connective. This is sound because
  implication is closed under disjunction, conjunction and quantification.
\item \emph{Symmetry}: the notion of derivation is top-down symmetric:
  a derivation is dualised essentially by flipping it upside-down.
  One example of how this symmetry is useful is the reduction of the
  cut rule to atomic form.
\end{enumerate}
The calculus of structures was conceived by Guglielmi in an earlier
unpublished version of \cite{Gug02}.  Its original purpose was to
express a logical system with a self-dual non-commutative connective
resembling sequential composition in process algebras
\cite{Gug02,GugStr01,GugStr02,PaolaBVL02}. The present work explores
the ideas developed in \cite{Gug02} in the setting of classical logic.
The calculus of structures has also been employed by Stra\ss burger in
\cite{StraELS01} to solve the problem of the non-local behaviour of
the promotion rule and in \cite{StraLSLL02} to give a local system for
full linear logic.  In the case of classical logic it led to a cut
elimination procedure similar to normalisation in natural deduction
\cite{BruACECL}.

This paper is structured as follows: in Section \ref{sec:struct-deriv}
I introduce the basic notions of the proof-theoretic formalism used, the
calculus of structures.  Section \ref{sec:propositional} is devoted to
classical propositional logic and Section \ref{sec:predicate} to
predicate logic.

The section for propositional logic is structured as follows: I first
present system $\SKSg$: a set of inference rules for classical
propositional logic, which is closed under a notion of duality.  I
translate derivations of a Gentzen-Schütte sequent system into this
system, and vice versa. This establishes soundness and completeness
with respect to classical propositional logic as well as cut
elimination. In the following I obtain an equivalent system, named
$\SKS$, in which identity, cut, weakening and contraction are reduced
to atomic form.  This entails locality of the system.  I go on to
establish three different normal forms for derivations, by what I call
`decomposition theorems'.

The outline of the section for predicate logic closely follows the one
for propositional logic. All results for the propositional systems
scale: reduction of cut, identity, weakening and contraction to atomic
form, cut elimination as well as the decomposition theorems.  The
resulting system with atomic rules is local except for the rules that
instantiate variables or check for free occurrences of a variable.

\section{The Calculus of Structures}\label{sec:struct-deriv}

\begin{Definition}
  \emph{Atoms} are denoted by $a$, $b$,~\dots.  The \emph{structures}
  of the language $\KS$ are generated by
  $$
  S \grammareq \false \mid \true \mid a \mid
  \pars{\,\underbrace{S,\dots,S}_{{}>0}\,} \mid
  \aprs{\,\underbrace{S,\dots,S}_{{}>0}\,} \mid \neg S \quadcm
  $$
  where $\false$ and $\true$ are the units \emph{false} and
  \emph{true}, ${\pars{S_1,\dots,S_h}}$ is a \emph{disjunction} and
  ${\aprs{S_1,\dots,S_h}}$ is a \emph{conjunction}. ${\bar S}$ is the
  \emph{negation} of the structure $S$. The negation of an atom is
  again an atom. Structures are denoted by $S$, $P$, $Q$, $R$, $T$,
  $U$, $V$ and $W$.  \emph{Structure contexts}, denoted by
  ${S\cons{\enspace}}$, are structures with one occurrence of
  $\cons{\enspace}$, the \emph{empty context} or \emph{hole}, that
  does not appear in the scope of a negation.  $S\cons R$ denotes the
  structure obtained by filling the hole in ${S\cons{\enspace}}$ with
  $R$. We drop the curly braces when they are redundant: for example,
  $S\pars{R,T}$ is short for $S\cons{\pars{R,T}}$.  A structure $R$ is
  a \emph{substructure} of a structure $T$ if there is a context
  $S\cons{\enspace}$ such that $S\cons{R}$ is $T$.
%Two substructures of a structure are called \emph{disjoint} 
%  if their sets of atom occurrences are disjoint.  
  Structures are \emph{equivalent} modulo the smallest equivalence
  relation induced by the equations shown in Figure~\ref{fig:Equations}.
  There, $\vec R, \vec T$ and $\vec U$ are finite sequences of
  structures, $\vec T$ is non-empty.  Structures are in \emph{normal
    form} if negation occurs only on atoms, and extra units as well as
  connectives are removed using the laws for units and associativity.
  In general we consider structures to be in normal form and do not
  distinguish between two equivalent structures.
\end{Definition}

\begin{Example}
  The structures $\overline{\pars{a,\false,b}}$ and $\aprs{\neg a,
    \true,\neg b}$ are equivalent, but they are not normal;
  $\aprs{\neg a,\neg b}$ is equivalent to them and normal, as well as
  $\aprs{\neg b, \neg a}$. The atom $a$ is not a substructure of
  $\overline{\pars{a,\false,b}}$, but $\neg a$ is a substructure of
  $\aprs{\neg a,\neg b}$.
\end{Example}

Structures are somewhere between formulae and sequents. They share
with formulae their tree-like shape and with sequents the built-in,
decidable equivalence modulo associativity and commutativity. The
equations for negation are adopted also in one-sided sequent systems,
so, apart from the equations for the units, the calculus of structures
does not use new equations. However, from the viewpoint of the
sequent calculus, it does extend the applicability of equations from
the level of sequents to the level of formulae.

\renewcommand{\arraystretch}{1.3}
\begin{figure}[t]
  \begin{center}
    \parbox[t]{\textwidth}{
      \parbox[t]{0.5\textwidth}{
        \smalltitle{Associativity}%
        $$
              \begin{array}[t]{c}
          \pars{\vec R,\pars{\vec T}, \vec U} =  
          \pars{\vec R,\vec T, \vec U}  \\
         \aprs{\vec R,\aprs{\vec T}, \vec U} =  
          \aprs{\vec R,\vec T, \vec U}  \\
        \end{array}
        $$
        \smalltitle{Units}%
        $$ 
        \begin{array}[t]{cc}
          \aprs{\false,\false} = \false \quad &
          \pars{\false, R} = R\\
          \pars{\true,\true} = \true  \quad &
          \aprs{\true, R} =  R\\
        \end{array}
        $$
        \smalltitle{Context Closure}%
%        $$ 
%        \begin{array}[h]{c}
%          S\cons{R}=S\cons{T}\\
%          \neg R = \neg T 
%        \end{array}
%          \quad \mbox{if} \quad R=T
%        $$
        $$ 
       \mbox{if} \quad R=T  \quad \mbox{then} \quad 
        \begin{array}[c]{c}
          S\cons{R}=S\cons{T} \\
          \neg R = \neg T
      \end{array}
        $$
        } 
      \parbox[t]{0.5\textwidth}{
        \smalltitle{Commutativity}%
        $$
        \begin{array}[t]{c}
          \pars{ R, T} = \pars{ T, R}\\
          \aprs{ R, T} = \aprs{ T, R}
        \end{array}
        $$
       \smalltitle{Negation}%
       $$ 
        \begin{array}[t]{c}
          \overline{\false}= \true \\ 
          \overline{\true}= \false \\      
          \overline{\strut\pars{R,T}} =
          \aprs{\neg R,\neg T}\\
           \overline{\strut\aprs{R,T}} =
          \pars{\neg R,\neg T}\\[1ex]
           {\skew3\neg{\neg R}= R} \\
        \end{array}
        $$
      } } 
    \caption{Syntactic equivalence of structures}
    \label{fig:Equations}
  \end{center}
\end{figure}

\begin{Definition}
  An \emph{inference rule\/} is a scheme of the kind
  $$
  \vcenter{\inf{\rho}
    {U}
    {V}}
  \quad,
  $$
  where $\rho$ is the \emph{name\/} of the rule, $V$ is its
  \emph{premise\/} and $U$ is its \emph{conclusion}.  If $V$ is of the
  form $S\cons{T}$ and $U$ is of the form $S\cons{R}$ then the
  inference rule is called \emph{deep}, otherwise it is called
  \emph{shallow}. In an instance of a deep inference rule
  $$
  \vcenter{\inf{\pi}
    {S\cons{R}}
    {S\cons{T}}}
  \quad,
  $$
  the structure taking the place of $R$ is its \emph{redex}, the
  structure taking the place of $T$ is its \emph{contractum} and the
  context taking the place of $S\cons{\enspace}$ is its
  \emph{context}.  A \emph{{\rm(}deductive\/{\rm)} system} $\sysS$ is
  a set of inference rules.
\end{Definition}

Most inference rules we will consider are deep. A deep inference rule
can be seen as a rewrite rule with the context made explicit. For
example, the rule $\pi$ from the previous definition seen top-down
corresponds to a rewrite rule $T\rightarrow R$. A shallow inference
rule can be seen as a rewrite rule that may only be applied to the
whole given term, not to arbitrary subterms.

\begin{Notation}
  To clarify the use of the syntactic equivalence where it is not
  obvious, I use the \emph{equivalence rule}
  $$\vcenter{\infnote{=}{R}{T}{\quad ,}}$$
  where $R$ and $T$ are equivalent
  structures.
\end{Notation}

\begin{Definition}
  A \emph{derivation} $\Delta$ in a certain deductive system
  is a finite sequence of instances of inference rules in the system:
  $$\vcenter{
    \der {\rho  }{R}  {
      \root{\rho' }{U                            } {
      \root{\pi' }{\vcenter{\hbox{\strut\vdots}}                            } {
        \root{\pi}{V                            }{
          \leaf        {T}}}}}
    }\quad.
    $$
    A derivation can consist of just one structure.  The topmost
    structure in a derivation is called the \emph{premise\/} of the
    derivation, and the structure at the bottom is called its
    \emph{conclusion}.  
\end{Definition}

Note that the notion of derivation is top-down symmetric, contrary to
the corresponding notion in the sequent calculus.

\begin{Notation}
  A derivation $\Delta$ whose premise is $T$, whose conclusion is $R$,
  and whose inference rules are in $\sysS$ is denoted by
  $$\vcenter{\xy\xygraph{[]!{0;<2pc,0pc>:}
      {T}-@{=}^<>(.5){\strut\sysS} _<>(.5){\strut\Delta}[d] {R}
    }\endxy}\quad .$$
\end{Notation}

\begin{Definition}
  A rule $\rho$ is \emph{derivable\/} for a system $\sysS$ if for
  every instance of $\vcenter{\inf{\rho}{R}{T}}\;$ there is a
  derivation $\vcenter{\xy\xygraph{[]!{0;<2pc,0pc>:}
      {T}-@{=}^<>(.5){\strut\sysS} _<>(.5){\strut}[d] {R} }\endxy}$.
\end{Definition}

The symmetry of derivations, where both premise and conclusion are
arbitrary structures, is broken in the notion of \emph{proof}:
\begin{Definition}
  A \emph{proof\/} is a derivation whose premise is the unit true. A
  proof $\Pi$ of $R$ in system $\sysS$ is denoted by
  $$\vcenter{\xy \xygraph{[]!{0;<2pc,0pc>:}
      {}*=<0pt>{}:@{|=}^<>(.5){\strut\sysS} _<>(.5){\Pi} [d] {R} }
    \endxy}\quad .$$
\end{Definition}

\section{Propositional Logic}\label{sec:propositional}
\subsection{A Symmetric System}\label{sec:sksg}

The following notion of duality is known as \emph{contrapositive}:

\begin{Definition}
  The \emph{dual} of an inference rule is obtained by exchanging
  premise and conclusion and replacing each connective by its De
  Morgan dual. For example
  $$
  \vcinf{\ruleidown} {S\pars{R,\neg R}} {S\cons{\true}} \qquad
  \mbox{is dual to} \qquad \vcinf{\ruleiup} {S\cons{\false}}
  {S\aprs{R,\neg R}} \quad ,
  $$
  where the rule $\ruleidown$ is called \emph{identity} and the
  rule $\ruleiup$ is called \emph{cut}.
\end{Definition}

The rules $\ruleidown$ and $\ruleiup$ respectively correspond to the
identity axiom and the cut rule in the sequent calculus, as we will
see shortly.

\begin{Definition}
  A system of inference rules is called \emph{symmetric} if for each
  of its rules it also contains the dual rule.
\end{Definition}

An example of a symmetric system is shown in Figure~\ref{fig:SKSg}. It
is called system $\SKSg$, where the first $\Sss$ stands for
``symmetric'', $\Kss$ stands for ``klassisch'' as in Gentzen's $\LK$
and the second $\Sss$ says that it is a system on structures. Small
letters are appended to the name of a system to denote variants. In
this case, the {\sf g} stands for ``general'', meaning that rules are
not restricted to atoms: they can be applied to arbitrary structures.
We will see in the next section that this system is sound and complete
for classical propositional logic.

\newcommand{\rulebox}[1]{\parbox{7em}{$$#1$$}}
\begin{figure}[htb]
  \begin{center}
    \fbox{
      \parbox{19em}{
        \rulebox{\vcinf{\ruleidown}
          {S\pars{R,\neg R}}
          {S\cons{\true}}}
        \hfill
        \rulebox{\vcinf{\ruleiup}
          {S\cons{\false}}
          {S\aprs{R,\neg R}}}
        $$
        \vcinf{\swir}
        {S\pars{\aprs{R,T},U }}
        {S\aprs{\pars{R,U},T}}
        $$
         \rulebox{\vcinf{\rulewdown}
          {S\cons{R}}
          {S\cons{\false}}}
        \hfill
        \rulebox{\vcinf{\rulewup}
          {S\cons{\true}}
          {S\cons{R}}}
        \rulebox{\vcinf{\rulecdown}
          {S\cons{ R}}
          {S\pars{ R,R }}}
        \hfill
        \rulebox{\vcinf{\rulecup}
          {S\aprs{ R,R }}
          {S\cons{ R}}}
        }
      }    
    \caption{System $\SKSg$}
    \label{fig:SKSg}
  \end{center}
\end{figure}

The rules $\swir, \rulewdown$ and $\rulecdown$ are called respectively
\emph{switch}, \emph{weakening} and \emph{contraction}.  Their dual
rules carry the same name prefixed with a ``co-'', so \eg $\rulewup$
is called \emph{co-weakening}.  Rules $\ruleidown, \rulewdown,
\rulecdown$ are called \emph{down-rules} and their duals are called
\emph{up-rules}. The dual of the switch rule is the switch rule
itself: it is \emph{self-dual}.

I now try to give an idea on how the familiar rules of the sequent
calculus correspond to the rules of $\SKSg$. For the sake of
simplicity I consider the rules of the sequent calculus in isolation,
\ie not as part of a proof tree. The full correspondence is shown in
Section \ref{sec:equiv}.

The identity axiom of the sequent calculus corresponds to the identity
rule $\ruleidown$: 
\renewcommand{\rulebox}[1]{\parbox{.33\textwidth}{$$#1$$}}
$$
\rulebox{\inf{}{\vdash  A, \neg A}{}}
\mbox{corresponds to}
\rulebox{\infnote{\ruleidown}{\pars{A,\neg A}}{\true}{\quad .}}
$$
However, $\ruleidown$ can appear anywhere in a proof, not only at
the top. The cut rule of the sequent calculus corresponds to the rule
$\ruleiup$ followed by two instances of the switch rule:
$$
\rulebox{\iinfnote{\sruleCut} 
             {\sqn{\Phi,\Psi}} {\sqn{\Phi, A}} 
             {\sqn{\Psi, \neg A}}{}}
\mbox{corresponds to}
\rulebox{\dernote{=}{}{\pars{\Phi,\Psi}}{
\root{\ruleiup}{\pars{\Phi,\Psi,\false}}{
\root{\rules}{\pars{\Phi,\Psi,\aprs{A,\neg A}}}{
\root{\rules}{\pars{\Phi,\aprs{A,\pars{\Psi,\neg
A}}}}{\leaf{\aprs{\pars{\Phi,A},\pars{\Psi,\neg A}}}}}
}}}.
$$
The multiplicative (or context-splitting) $\sruleRand$ in the
sequent calculus corresponds to two instances of switch rule:
$$
\rulebox{\iinfnote{\sruleRand} {\sqn{\Phi,\Psi,A \wedge B}}
  {\sqn{\Phi, A}} {\sqn{\Psi, B}}{}}
\mbox{corresponds to}
\rulebox{\dernote{\rules}{}{\pars{\Phi,\Psi,\aprs{A,B}}}{
    \root{\rules}{\pars{\Phi,\aprs{A,\pars{\Psi,B}}}}{\leaf{\aprs{\pars{\Phi,A},\pars{\Psi,B}}}
    }}}.
$$

A contraction in the sequent calculus corresponds to the
$\rulecdown$ rule:
$$
\rulebox{\infnote{\sruleRC} 
             {\sqn{\Phi,A}} 
             {\sqn{\Phi,A,A}}{}}
\mbox{corresponds to}
\rulebox{\infnote{\rulecdown}
          {\pars{\Phi,A}}
          {\pars{\Phi,A,A}}{\quad ,}}
$$

just as the weakening in the sequent calculus corresponds
to the $\rulewdown$ rule: 
$$
\rulebox{\infnote{\sruleRW} 
             {\sqn{\Phi,A}} 
             {\sqn{\Phi}}{}}
\mbox{corresponds to}
\rulebox{\dernote{\rulewdown}{\quad .}
          {\pars{\Phi,A}}
          {\root{=}{\pars{\Phi,\false}}{\leaf{\Phi}}}}
$$

The $\rulecup$ and $\rulewup$ rules have no analogue in the sequent
calculus. Their role is to ensure that our system is symmetric.  They
are obviously sound since they are just duals of the rules
$\rulecdown$ and $\rulewdown$ which correspond to sequent calculus
rules.

Derivations in a symmetric system can be dualised:

\begin{Definition}
  The \emph{dual} of a derivation is obtained by turning it
  upside-down and replacing each rule, each connective and each atom
  by its dual.
\end{Definition}

For example
$$\downsmash{ \dernote{\rulecdown}{\qquad\quad \mbox{is dual to}\quad
    \qquad}{a} { \root {\rulewup}
    {\pars{ a,a }} { \leaf {\pars{\aprs{a, \neg b}, a}} }}
  \dernote{\rulewdown}{\qquad.} {\aprs{\pars{\neg a, b},\neg a}} {
    \root{\rulecup}{\aprs{\neg a,\neg a}} { \leaf{\neg a}}}}
$$

The notion of proof, however, is an asymmetric one: the dual of a proof is
not a proof.

\hyphenation{Sequent}
\hyphenation{Calculus}
\subsection{Correspondence to the Sequent Calculus}
\label{sec:equiv}

The sequent system that is most similar to system $\SKSg$ is the
one-sided system $\GSonep$ \cite{TroSch96}, also called
\emph{Gentzen-Schütte} system. In this section we consider a version
of $\GSonep$ with multiplicative context treatment and constants
$\top$ and $\bot$, and we translate its derivations to derivations in
$\SKSg$ and vice versa. Translating from the sequent calculus to the
calculus of structures is straightforward, in particular, no new cuts
are introduced in the process. But to translate in the other direction
we have to simulate deep inferences in the sequent calculus, which is
done by using the cut rule.

One consequence of those translations is that system $\SKSg$ is sound
and complete for classical propositional logic. Another consequence is
cut elimination: one can translate a proof with cuts in $\SKSg$ to a
proof in $\GSonep+\sruleCut$, apply cut elimination for $\GSonep$, and
translate back the resulting cut-free proof to obtain a cut-free proof
in $\SKSg$.

\begin{Definition}
  System $\GSonep$ is the set of rules shown in Figure
  \ref{fig:GSOne}. The system $\GSonep+\sruleCut$ is $\GSonep$
together with
$$
      \vciinf{\sruleCut} 
             {\sqn{\Phi,\Psi}} {\sqn{\Phi, A}} 
             {\sqn{\Psi, \neg A}} \quad .
$$

\emph{Formulae} are denoted by $A$ and $B$. They contain negation only
on atoms and may contain the constants $\top$ and $\bot$. Multisets of
formulae are denoted by $\Phi$ and $\Psi$. The empty multiset is
denoted by $\emptyset$. In $A_1,\dots,A_h$, where
$h\geq 0$, a formula denotes the corresponding singleton multiset and
the comma denotes multiset union. \emph{Sequents}, denoted by
$\Sigma$, are multisets of formulae.  \emph{Derivations} are denoted
by $\Delta$ or
$\vcenter{\Derivation{\Sigma_1}{\enspace\cdots\enspace}{\Sigma_h}\Sigma\Delta}$,
where $h \geq 0$, the sequents $\Sigma_1,\dots,\Sigma_h$ are the
\emph{premises} and $\Sigma$ is the \emph{conclusion}. Proofs, denoted
by $\Pi$, are derivations where each leaf is an instance of
$\sruleAx$.
\end{Definition}

\renewcommand{\rulebox}[1]{\parbox{.29\textwidth}{$$#1$$}}
\begin{figure}[tb]
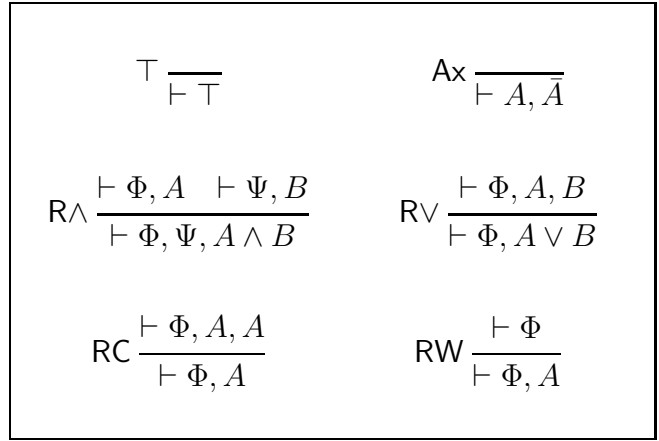

  \begin{center}
    \fbox{ \parbox{.6\textwidth}{\medskip
        \rulebox{\inf{\top}{\sqn{\top}}{}}
        \rulebox{\inf{\sruleAx}{\sqn{A,\neg A}}{}} 
\medskip
        \rulebox{\iinf{\sruleRand} {\sqn{\Phi, \Psi, A\wedge B}}
          {\sqn{\Phi, A}} {\sqn{\Psi, B}}} \rulebox{\inf{\sruleRv}
          {\sqn{\Phi, A\vee B}} {\sqn{\Phi, A, B}}} 
\medskip
        \rulebox{\inf{\sruleRC} {\sqn{\Phi, A}} {\sqn{\Phi, A, A}}}
        \rulebox{\inf{\sruleRW} {\sqn{\Phi, A}} {\sqn{\Phi}}} }}
    \caption{$\GSonep$: classical logic in Gentzen-Schütte form}
    \label{fig:GSOne}
  \end{center}
\end{figure}

\subsubsection*{From the Sequent Calculus to the Calculus of Structures}

\begin{Definition}
  The function $\ls{\enspace . \enspace}$ maps formulae, multisets of
  formulae and sequents of $\GSonep$ to structures:
$$
\begin{array}{cccl}
\ls a  & = & a\;,      \\[1ex] 
\ls \top  & = & \true\;,      \\[1ex] 
\ls \bot  & = & \false\;,      \\[1ex] 
\ls {A\lor B} & = & \pars{\ls A, \ls B}\;, \\[1ex]
\ls {A\land B} & = & \aprs{\ls A, \ls B}\;, \\[2ex] 
\ls{\emptyset} & = & \false\;, \\[1ex] 
\ls {A_1,\dots,A_h} & = & \pars{\ls {A_1}, \dots, \ls {A_h}}\;, & 
%\ls {\sqn{A_1,\dots,A_h}} & = & \pars{\ls {A_1}, \dots, \ls {A_h}}, &
\quad \mbox{where } h>0.\\
\end{array}
$$
\end{Definition}

In proofs, when no confusion is possible, the subscript \ls{\enspace}
may be dropped to improve readability.  

In the following, we will put derivations of system $\SKSg$ into a
context. This is possible because all inference rules of the system
are deep.

\begin{Definition}
  Given a derivation $\Delta$, the derivation $S\cons{\Delta}$ is
obtained as follows:
  $$\Delta = \vcenter{ \der {\rho }{R} { \root{\rho' }{U } { \root{\pi
        }{\vcenter{\hbox{\strut\vdots}} } { \root{\pi'}{V }{ \leaf
            {T}}}}} }\qquad\qquad S\cons{\Delta} = \vcenter{ \der
    {\rho }{S\cons{R}} { \root{\rho' }{S\cons{U} } { \root{\pi
        }{\vcenter{\hbox{\strut\vdots}} } { \root{\pi'}{S\cons{V} }{
            \leaf {S\cons{T}}}}}} }\qquad .
  $$
\end{Definition}

\newcommand{\KSgcut}{\SKSg\,\setminus\,\{\rulecup,\rulewup\}}
\medskip
\begin{Theorem}\label{thm:GS1toSKS}
  For every derivation $
  \vcenter{\Derivation{\Sigma_1}{\enspace\cdots\enspace}{\Sigma_h}\Sigma{}}
  \; \mbox{in $\GSonep+\sruleCut$} $ there exists a derivation
  $
  \vcenter{ \xy \xygraph{[]!{0;<3pc,0pc>:}
      {\aprs{{\ls{\Sigma_1}},\dots,\ls{\Sigma_h}}}-@{=}^<>(.5){\strut\KSgcut}[d]
      { \ls{\Sigma} } } \endxy }$ 
 with the same number of cuts.
\end{Theorem}
\begin{proof}
  By structural induction on the given derivation $\Delta$.
  \\[2ex]
\noindent \smalltitle{Base Cases}\\
\begin{enumerate}
\item $ \Delta = \Sigma\,.$ Take $\ls{\Sigma}$.
\item $\Delta = \vcinf{\top}{\sqn{\top}}{}\,.$ Take
  $\true\,.$
\item $\Delta = \vcinf{\sruleAx}{\sqn{A,\neg A}}{}\,.$ Take
  $\vcinf{\intr}{\pars{ \ls{A}, {\neg {\ls{A}}}}}{\true}\,.$
\end{enumerate}

\noindent \smalltitle{Inductive Cases}

In the case of the $\sruleRand$ rule, we have a derivation 
$$\Delta \; = \; \vcenter{ \ddernote{\sruleRand}{}{\sqn{\Phi,\Psi, A\land B}}{
    \Derivationleaf {\Sigma_1} {\enspace\cdots\enspace} {\Sigma_k}
    {\sqn{\Phi, A}} {} } { \Derivationleaf {\Sigma'_1}
    {\enspace\cdots\enspace} {\Sigma'_l} {\sqn{\Psi, B}} {} } } \quad .
$$  

By induction hypothesis we obtain derivations 

$$
\vcenter{\xy \xygraph{[]!{0;<3pc,0pc>:}
  {\aprs{{{\Sigma_1}},\dots,{\Sigma_k}}}-@{=}_{\Delta_1}^<>(.5){\strut\KSgcut}[d]
  { \pars{{\Phi},{A}} } } \endxy}
\qquad \mbox{and} \qquad
\vcenter{\xy \xygraph{[]!{0;<3pc,0pc>:}
  {\aprs{{{\Sigma'_1}},\dots,{\Sigma'_l}}}-@{=}_{\Delta_2}^<>(.5){\strut\KSgcut}[d]
  { \pars{{\Psi},{B}} } } \endxy} \quad .
$$
The derivation $\Delta_1$ is put into the context
$\aprs{\cons{\enspace},{\Sigma'_1},\dots,{\Sigma'_l}}$ to obtain
$\Delta'_1$ and the derivation $\Delta_2$ is put into the context
$\aprs{\pars{{\Phi},{A}},\cons{\enspace}}$ to obtain $\Delta'_2$:
$$
\vcenter{\xy \xygraph{[]!{0;<3pc,0pc>:}
  {\aprs{{\Sigma_1},\dots,{\Sigma_k},{\Sigma'_1},\dots,{\Sigma'_l}}}-@{=}_{\Delta'_1}^<>(.5){\strut\KSgcut}[d]
  { \aprs{\pars{{\Phi},{A}},{\Sigma'_1},\dots,{\Sigma'_l}} } } \endxy}
\qquad \mbox{and} \qquad
\vcenter{\xy \xygraph{[]!{0;<3pc,0pc>:}
  {\aprs{\pars{{\Phi},{A}},{\Sigma'_1},\dots,{\Sigma'_l}}}-@{=}_{\Delta'_2}^<>(.5){\strut\KSgcut}[d]
  { \aprs{\pars{{\Phi},{A}},\pars{{\Psi},{B}}} } } \endxy} \quad .
$$
The derivation in $\SKSg$ we are looking for is obtained by composing
$\Delta'_1$ and $\Delta'_2$ and applying the switch rule twice:

\renewcommand{\rulebox}[1]{\parbox{\textwidth}{$$#1$$}}

$$
\vcenter{\xy\xygraph{[]!{0;<1pc,0pc>:}
    {\aprs{{\Sigma_1},\dots,{\Sigma_k},{\Sigma'_1},\dots,{\Sigma'_l}}}
    -@{=}_{\Delta'_1}^<>(.5){\strut\KSgcut}[dddd]
    {\aprs{\pars{{\Phi},{A}},{\Sigma'_1},\dots,{\Sigma'_l}}}
    -@{=}_<>(.5){\Delta'_2}^<>(.5){\strut\KSgcut}[dddddd] 
      { \dernote
      {\rules }{\quad .}{\pars{{\Phi},{\Psi},\aprs{{A},{B}}}} {
        \rootnote{\quad \rules }{}{
          \pars{{\Psi},\aprs{\pars{{\Phi},{A}},{B}}} } {
          \leaf {
            \aprs{\pars{{\Phi},{A}},\pars{{\Psi},{B}}} }
        }} } }\endxy 
}
$$

The other cases are similar. The only case that requires a cut in
$\SKSg$ is a cut in $\GSonep$.
\end{proof}

\begin{Corollary}\label{cor:gssks}
  If a sequent $\Sigma$ has a proof in $\GSonep+\sruleCut$ then $\ls{\Sigma}$
  has a proof in $\KSgcut$.
\end{Corollary}

\begin{Corollary}\label{cor:gssks2}
  If a sequent $\Sigma$ has a proof in $\GSonep$ then
  $\ls{\Sigma}$ has a proof in
  $\SKSg\,\setminus\,\{\ruleiup,\rulecup,\rulewup\}$.
\end{Corollary}

\subsubsection*{From the Calculus of Structures to the Sequent Calculus}

In the following, structures are assumed to contain 1) negation only
on atoms and 2) only conjunctions and disjunctions of exactly two
structures. That is not a restriction because for each structure there
exists an equivalent one which has these properties.

\begin{Definition}
  The function $\lg{\enspace . \enspace}$ maps structures of $\SKSg$ to
  formulae of $\GSonep$:
$$
\matrix{
    \lg a & = & a,     \quad \cr\vspace{\jot}
    \lg \true & = & \top,   \quad \cr\vspace{\jot}
    \lg \false & = & \bot, \quad \cr\vspace{\jot}
    \lg {\pars{R, T}} & = & \lg R \lor \lg T, \cr\vspace{\jot}
    \lg {\aprs{R, T}} & = & \lg R \land \lg T. \cr
}
$$
\end{Definition}

\begin{Lemma}\label{lem:context}
  For every two formulae $A,B$ and every formula context
  $C\cons{\enspace}$ there exists a derivation
  $\vcenter{\Derivation{\quad}{\sqn{A,\neg
        B}}{\quad}{\sqn{C\cons{A},\overline{C\cons{B}}}}{}}$ in
  $\GSonep$.
\end{Lemma}

\begin{proof}
  By structural induction on the context $C\cons{\enspace}$. The base
  case in which $C\cons{\enspace}=\cons{\enspace}$ is trivial. If
  $C\cons{\enspace}=C_1\wedge C_2\cons{\enspace}$, then the derivation
  we are looking for is
  $$
  \dernote{\sruleRv}{\quad ,}{\sqn{C_1 \wedge C_2\cons{A},\neg C_1
      \vee \overline{C_2\cons{B}}}}{\rroot{\sruleRand}{\sqn{C_1 \wedge
        C_2\cons{A},\neg C_1,
        \overline{C_2\cons{B}}}}{\root{\sruleAx}{\sqn{C_1,\neg
          C_1}}{\leaf{}}}{\leaf{\Derivation{\quad}{\sqn{A,\neg
            B}}{\quad}{\sqn{C_2\cons{A},\overline{C_2\cons{B}}}}\Delta
      }}}
  $$
  where $\Delta$ exists by induction hypothesis. The other case, in
  which $C\cons{\enspace}=C_1\vee C_2\cons{\enspace}$, is similar.
\end{proof}

\begin{Theorem}\label{thm:SKStoGS1} For every derivation $\vcenter{ \xy
    \xygraph{[]!{0;<3pc,0pc>:} {Q}-@{=}^<>(.5){\strut\SKSg}[d] {P} }
    \endxy }$ there exists a derivation
  $\vcenter{\Derivation{\quad\quad}{\sqn{\lg{Q}}}{\quad\quad}{\sqn{\lg{P}}}{\,}}$
  in $\GSonep+\sruleCut$.
\end{Theorem}

\begin{proof} We construct the sequent derivation by
  induction on the length of the given derivation $\Delta$ in $\SKSg$.\\[1ex]
\noindent \smalltitle{Base Case} 

If $\Delta$ consists of just one structure $P$, then $P$ and $Q$ are
the same. Take $\sqn{\lg{P}}$.
\\[2ex]
\noindent \smalltitle{Inductive Cases}\nopagebreak

We single out the topmost rule instance in $\Delta$:
$$
\vcenter{ \xy \xygraph{[]!{0;<3pc,0pc>:}
    {Q}-@{=}_\Delta^<>(.5){\strut\SKSg}[d] {P} } \endxy } = \vcenter{
  \xy \xygraph{[]!{0;<2pc,0pc>:}
    {\inf{\rho}{S\cons{R}}{S\cons{T}}}-@{=}_<>(.5){\Delta'}^<>(.5){\strut\SKSg}[dd]
    {P} } \endxy }
$$

The corresponding derivation in $\GSonep$ will be as follows:

$$
\vcenter{
\iinfnote{\sruleCut}{\Derivation{}{\qquad\sqn{S\cons{R}}\qquad}{}{\sqn{P}}{\Delta_2}}{
\Derivation{\quad}{ 
\Derivation{\quad}{\hspace{10ex}}{\quad}{\sqn{R,\neg T}}{\Pi}
}
{}{\sqn{S\cons{R},\overline{S\cons{T}}}}{\Delta_1}
}{\sqn{S\cons{T}}}{\quad ,} 
}
$$

where $\Delta_1$ exists by Lemma \ref{lem:context} and $\Delta_2$
exists by induction hypothesis. The proof $\Pi$ depends on the rule
$\rho$. In the following we will see that the proof $\Pi$ exists for all
the rules of $\SKSg$.

For identity and cut, \ie
$$
\vcinf{\ruleidown}
          {S\pars{U,\neg U}}
          {S\cons{\true}}
        \qquad \mbox{and} \qquad
\vcinf{\ruleiup}
          {S\cons{\false}}
          {S\aprs{U,\neg U}} \quad ,
$$
we have the following proofs:
$$
\vcenter{
\dernote{\sruleRW}{}
        {\sqn{U \vee \neg U, \bot}}
             {\root{\sruleRv}
                    { \sqn{U \vee \neg U} }
                     {\root{\sruleAx}
                          {\sqn{ U,\neg U } 
                          }
                          {\leaf{}
                          }
                    } 
             }  
        }
\qquad \mbox{and} \qquad
\vcenter{
\dernote{\sruleRW}{}
        {\sqn{\bot,\neg U \vee U}}
             {\root{\sruleRv}
                    { \sqn{\neg U \vee U} }
                     {\root{\sruleAx}
                          {\sqn{ U,\neg U } 
                          }
                          {\leaf{}
                          }
                    } 
             }  
        } \quad .
$$

In the case of the switch rule, \ie 
$$
\vcinf{\swir} {S\pars{\aprs{U,T},V }} {S\aprs{\pars{U,V},T}} \quad
,
$$
we have
$$
\vcenter{ \dernote{\sruleRv^2}{} {\sqn{ (U \wedge T) \vee V, (\neg
      U \wedge \neg V) \vee \neg T }} {\rroot{\sruleRand} { \sqn{ (U
        \wedge T) , V, (\neg U \wedge \neg V) , \neg T } }
    {\rroot{\sruleRand} {\sqn{ U,\neg U \wedge \neg V, V } }
      {\root{\sruleAx} {\sqn{U,\neg U } } {\leaf{} } }
      {\root{\sruleAx} {\sqn{ V, \neg V } } {\leaf{} } } }
    {\root{\sruleAx} {\sqn{ T,\neg T } } {\leaf{} } } } } \quad .
$$  

For contraction and its dual, \ie
$$
\vcinf{\rulecdown} {S\cons{ U}} {S\pars{ U,U }} \qquad \mbox{and}
\qquad \vcinf{\rulecup} {S\aprs{ U,U }} {S\cons{ U}} \quad ,
$$
we have 
$$
\vcenter{ \dernote{\sruleRC}{} {\sqn{U,\neg U \wedge \neg U}}
  {\rroot{\sruleRand} { \sqn{U,U,\neg U \wedge \neg U} }
    {\root{\sruleAx} {\sqn{U,\neg U } } {\leaf{} } } {\root{\sruleAx}
      {\sqn{ U,\neg U } } {\leaf{} } } } } \qquad \mbox{and} \qquad
\vcenter{ \dernote{\sruleRC}{} {\sqn{U \wedge U, \neg U}}
  {\rroot{\sruleRand} { \sqn{U \wedge U,\neg U, \neg U} }
    {\root{\sruleAx} {\sqn{U,\neg U } } {\leaf{} } } {\root{\sruleAx}
      {\sqn{ U,\neg U } } {\leaf{} } } } } \quad .
$$
For weakening and its dual, \ie
$$
\vcinf{\rulewdown} {S\cons{U}} {S\cons{\false}} \qquad \mbox{and}
\qquad \vcinf{\rulewup} {S\cons{\true}} {S\cons{U}}\quad ,
$$
we have 
$$
\vcenter{ \dernote{\sruleRW}{} {\sqn{U, \top}} {\root{\top} {
      \sqn{\top} } {\leaf{} } } } \qquad \mbox{and} \qquad \vcenter{
  \dernote{\sruleRW}{} {\sqn{\top, \neg U}} {\root{\top} { \sqn{\top}
    } {\leaf{} } } } \quad .
$$  
\end{proof}

\begin{Corollary}\label{cor:sksgs}
  If a structure $S$ has a proof in $\SKSg$ then $\sqn{\lg{S}}$ has a
  proof in $\GSonep+\sruleCut$.
\end{Corollary}

Soundness and completeness of $\SKSg$, \ie the fact that a structure
has a proof if and only if it is valid, follows from soundness and
completeness of $\GSonep$ by Corollaries \ref{cor:gssks} and
\ref{cor:sksgs}. Moreover, a structure $T$ implies a structure $R$ if
and only if there is a derivation from $T$ to $R$, which follows from
soundness and completeness and the following theorem.

\begin{Theorem}\label{thm:ded}\theoremnl
  There is a derivation $\vcenter{\xy\xygraph{[]!{0;<1pc,0pc>:} {T}
      -@{=}^<>(.5){\strut\SKSg}[ddd] {R} }\endxy}$ if and only
  if there is a proof $\vcenter{\xy\xygraph{[]!{0;<1pc,0pc>:} {}*=<0pt>{}:
      @{|=}^<>(.5){\strut\SKSg}[ddd] {\pars{\neg T,R}}
    }\endxy}\quad .$
\end{Theorem}

\begin{proof}
  A proof $\Pi$ can be obtained from a given derivation $\Delta$ and a
  derivation $\Delta$ from a given proof $\Pi$, respectively, as follows:
  $$\vcenter{\xy\xygraph{[]!{0;<1pc,0pc>:}
      {\inf{\ruleidown}{\pars{\neg T, T}}{\true}\hspace{2.7ex}}
      -@{=}^<>(.5){\strut\SKSg} _<>(.5){\strut\pars{\neg
          T,\Delta}}[dddd] {\pars{ \neg T, R}} }\endxy} \qquad
  \mbox{and} \qquad \vcenter{\xy\xygraph{[]!{0;<1pc,0pc>:} {T}
      -@{=}_<>(.5){\aprs{T,\Pi}}^<>(.5){\strut\SKSg} [ddddd]
      {\dernote{\aintr}{}{R}{ \root{\swir}{\pars{R,\aprs{T,\neg T}}}{
            \leaf{\aprs{T,\pars{\neg T,R}} }}}} }\endxy} \qquad .
$$
\end{proof}

\subsection{Admissibility of the Cut and the Other Up-Rules}
\label{sec:cutelim}

If one is just interested in provability, then the up-rules of the
symmetric system $\SKSg$, \ie $\ruleiup$, $\rulewup$ and $\rulecup$,
are superfluous.  By removing them we obtain the asymmetric, cut-free
system shown in Figure~\ref{fig:KSg}, which is called system $\KSg$.

\begin{Definition}
  A rule $\rho$ is \emph{admissible} for a system $\sysS$ if for every
  proof $\vcenter{\xy \xygraph{[]!{0;<2pc,0pc>:}
      {}*=<0pt>{}:@{|=}^<>(.5){\sysS\cup\{\rho\}} _<>(.5){} [d] {S} }
    \endxy} $ there is a proof $\vcenter{\xy
    \xygraph{[]!{0;<2pc,0pc>:} {}*=<0pt>{}:@{|=}^<>(.5){\sysS}
      _<>(.5){} [d] {S} } \endxy} $.
\end{Definition}

The admissibility of all the up-rules for system $\KSg$ is shown
by using the translation functions from the previous section:

\newcommand{\SKSproof}{\xy\xygraph{[]!{0;<1pc,0pc>:}
    {}*=<0pt>{}:@{|=}^<>(.5){\strut\SKSg}[dd] {S} }\endxy}
\newcommand{\KSproof}{\xy\xygraph{[]!{0;<1pc,0pc>:}
    {}*=<0pt>{}:@{|=}^<>(.5){\strut\KSg}[dd] {S} }\endxy}
\newcommand{\GScutproof}{
  \Derivation{\quad}{\hspace{9ex}}{\quad}{\sqn{\lg{S}}}{{\GSonep \atop
+\sruleCut}}}
\newcommand{\GSproof}{
  \Derivation{\quad}{\hspace{9ex}}{\quad}{\sqn{\lg{S}}}{\GSonep}}

\begin{Theorem}\label{thm:transcutelim}
  The rules $\ruleiup$, $\rulewup$ and $\rulecup$ are admissible for
  system $\KSg$.
\end{Theorem}

\begin{proof}
\theoremnl
\xymatrix@C11ex{
{\SKSproof} \ar[r]^-{\mbox{\scriptsize Corollary \ref{cor:sksgs}}} &  
{\GScutproof} \ar[r]^-{\mbox{\scriptsize Cut
elimination} \atop \mbox{ \scriptsize for $\GSonep$}} & 
{\GSproof} \ar[r]^-{\mbox{\scriptsize Corollary \ref{cor:gssks2} }} & 
{\KSproof} \\
}
\end{proof}

\begin{Definition}
  Two systems $\sysS$ and $\sysS'$ are (\emph{weakly})
  \emph{equivalent} if for every proof
  $\vcenter{\xy\xygraph{[]!{0;<2pc,0pc>:}
      {}*=<0pt>{}:@{|=}^<>(.5){\strut\sysS}[d] {R} }\endxy}$ there is a proof
  $\vcenter{\xy\xygraph{[]!{0;<2pc,0pc>:}
      {}*=<0pt>{}:@{|=}^<>(.5){\strut\sysS'}[d] {R} }\endxy}$, and vice versa.
\end{Definition}

\begin{Corollary}
  The systems $\SKSg$ and $\KSg$ are equivalent.
\end{Corollary}

\begin{Definition}
   Two systems
  $\sysS$ and $\sysS'$ are \emph{strongly equivalent} if for every
  derivation $\vcenter{\xy\xygraph{[]!{0;<2pc,0pc>:}
      {T}-@{=}^<>(.5){\strut\sysS}[d] {R} }\endxy}$ there is a
  derivation $\vcenter{\xy\xygraph{[]!{0;<2pc,0pc>:}
      {T}-@{=}^<>(.5){\strut\sysS'}[d] {R} }\endxy}$, and vice versa.
\end{Definition}

\begin{Remark}
  The systems $\SKSg$ and $\KSg$ are not strongly equivalent.  The cut
  rule, for example, can not be derived in system $\KSg$.
\end{Remark}

When a structure $R$ implies a structure $T$ then there is not
necessarily a derivation from $R$ to $T$ in $\KSg$, while there is one
in $\SKSg$. Therefore, I will in general use the asymmetric, cut-free
system for deriving conclusions from the unit true, while I will use
the symmetric system (\ie the system with cut) for deriving
conclusions from arbitrary premises.

\renewcommand{\rulebox}[1]{\mbox{$#1$}}
\begin{figure}[tb]
  \begin{center}
    \fbox{
      \parbox{\textwidth}{
$$
        \rulebox{\vcinf{\ruleidown}
          {S\pars{R,\neg R}}
          {S\cons{\true}}}
        \qquad\qquad
        \rulebox{\vcinf{\swir}
        {S\pars{\aprs{R,U},T }}
        {S\aprs{\pars{R,T},U}}}
        \qquad\qquad
        \rulebox{\vcinf{\rulewdown}
          {S\cons{R}}
          {S\cons{\false}}}
        \qquad\qquad
        \rulebox{\vcinf{\rulecdown}
          {S\cons{ R}}
          {S\pars{ R,R }}}
$$
        }
      }    
    \caption{System $\KSg$}
    \label{fig:KSg}
  \end{center}
\end{figure}

As a result of cut elimination, sequent systems fulfill the subformula
property.  Our case is different, because our rules do not split the
derivation according to the main connective of the active formula.
However, seen bottom-up, in system $\KSg$ no rule introduces new
atoms.  It thus satisfies one main aspect of the subformula
property: when given a conclusion of a rule there is only a finite
number of premises to choose from. In proof search, for example, the
branching of the search tree is finite.

There is also a semantic cut elimination proof for system $\SKSg$,
analogous to the one given in \cite{TroSch96} for system $\Gthree$.
The given proof with cuts is thrown away, keeping only the information
that its conclusion is valid, and a cut-free proof is constructed from
scratch. This actually gives us more than just admissibility
of the up-rules: it also yields a separation of proofs into distinct
phases.

\begin{Theorem}\label{thm:semcutelim}
\theoremnl
     For every proof $\vcenter{\xy\xygraph{[]!{0;<1pc,0pc>:}
         {}*=<0pt>{}:@{|=}^<>(.5){\strut\SKSg}[dd] {S} }\endxy}$ there
     is a proof $\vcenter{\xy\xygraph{[]!{0;<1pc,0pc>:}
         {}*=<0pt>{}:@{|=}^<>(.5){\strut\{\ruleidown\}}[dd] {S''}
         -@{=}^<>(.5){\strut\{\rulewdown\}}[dd] {S'}
         -@{=}^<>(.5){\strut\{\rules,\rulecdown\}} [dd] {S}
       }\endxy}$.
\end{Theorem}

\begin{proof}

Consider the rule \emph{distribute}:
$$
\infnote{\ruleddown}
        {S\pars{R,\aprs{T,U}}}
        {S\aprs{\pars{R,T},\pars{R,U}}}{\quad ,}
$$
\\
 which can be realized by a contraction and two switches:
 $$\downsmash{ \dernote {\rulecdown }{\quad .}{S\pars{R,\aprs{T,U}}} {
     \root{\rules}{S\pars{R,R,\aprs{T,U}}}{
       \root{\rules}{S\pars{R,\aprs{\pars{R,T},U}}}{
         \leaf{S\aprs{\pars{R,T},\pars{R,U}}}}}} }$$
 Build a
 derivation $\vcenter{\xy\xygraph{[]!{0;<1pc,0pc>:} {S'}
     -@{=}^<>(.5){\strut\{\ruleddown\}}[dd] {S} }\endxy}$, by going
upwards from $S$ applying $\ruleddown$ as many times as possible. Then
$S'$ will be in conjunctive normal form, \ie 

$$
S' = \aprs{\pars{a_{11}, a_{12},\dots}, \pars{a_{21}, a_{22},\dots
    }, \dots, \pars{a_{n1}, a_{n2},\dots}}\quad.
$$
$S$ is valid because there is a proof of it. The rule $\ruleddown$ is
invertible, so $S'$ is also valid. A conjunction is valid only if all
its immediate substructures are valid. Those are disjunctions of
atoms. A disjunction of atoms is valid only if it contains an atom
$a$ together with its negation $\neg a$. Thus, more specifically, $S'$ is of the
form
$$
S' = \aprs{\pars{b_1,\neg b_1,a_{11}, a_{12},\dots}, \pars{b_2,\neg
    b_2,a_{21}, a_{22},\dots }, \dots, \pars{b_n,\neg b_n,a_{n1},
    a_{n2},\dots}}\quad .
$$
Let $ S'' = \aprs{\pars{b_1,\neg b_1}, \pars{b_2,\neg b_2}, \dots,
  \pars{b_n,\neg b_n}} \quad$.

Obviously, there is a derivation
$\vcenter{\xy\xygraph{[]!{0;<1pc,0pc>:} {S''}
    -@{=}^<>(.5){\strut\{\rulewdown\}}[dd] {S'} }\endxy}$ and a proof
$\vcenter{\xy\xygraph{[]!{0;<1pc,0pc>:}
    {}*=<0pt>{}:@{|=}^<>(.5){\strut\{\ruleidown\}}[dd] {S''}
    }\endxy}$.  
\end{proof}

\subsection{Reducing Rules to Atomic Form}\label{sec:atomic}

In the sequent calculus, the identity rule can be reduced to its
atomic form. The same is true for our system, \ie
$$
\vcinf{\intr} {S\pars{R,\neg R}} {S\cons{\true}} \qquad \mbox{is
equivalently replaced by} \qquad
\vcinf{\ruleaidown} {S\pars{a,\neg a}} {S\cons{\true}} \quad ,
$$
where $\ruleaidown$ is the \emph{atomic identity} rule. Similarly
to the sequent calculus, this is achieved by inductively
replacing an instance of the general identity rule by instances on
smaller structures:
$$
\vcinf{\intr} {S\pars{P,Q,\aprs{\neg P,\neg Q}}} {S\cons{\true}} \qquad \leadsto \qquad
\vcenter{
\dernote{\swir}{}{S\pars{P,Q,\aprs{        \neg P ,        \neg Q }}} {
\root   {\swir}  {S\pars{  Q,\aprs{\pars{P,\neg P},        \neg Q }}} {
\root   {\intr}  {S    \aprs{\pars{P,\neg P},\pars{Q,\neg Q}} }{
\root   {\intr}  {S        \pars{Q,\neg Q} }{
\leaf            {S\cons{\true}}                                   }}}}
}
\qquad .
$$

What is new in the calculus of structures is that the cut can also be
reduced to atomic form: just take the dual derivation of
the one above: 

$$
\vcinf{\ruleiup}{S\cons{\false}}{S\aprs{\neg P,\neg Q,\pars{P,Q}}}
\qquad \leadsto \qquad \vcenter{ \dernote{\ruleiup}{}{S\cons{\false}}
  { \root {\ruleiup} {S\aprs{\neg Q, Q}} { \root {\rules}
      {S\pars{\aprs{\neg P, P},\aprs{\neg Q, Q}} }{ \root {\rules}
        {S\aprs{\neg Q,\pars{\aprs{\neg P, P},Q}} }{ \leaf
          {S\aprs{\neg Q,\neg P,\pars{Q, P}} } }}}}} \qquad .
$$

It turns out that weakening can also be reduced to atomic form. When
identity, cut and weakening are restricted to atomic form, there is only
one non-local rule left in system $\KSg$: contraction. It can not be
reduced to atomic form in system $\KSg$. Tiu solved 
this problem in when he discovered the \emph{medial} rule \cite{BruTiu01}:
$$
\vcinf{\rulem} {S\aprs{\pars{R,T},\pars{U,V} }}
{S\pars{\aprs{R,U},\aprs{T,V} }} \quad .
$$

This rule has no analogue in the sequent calculus. But it is
clearly sound because we can derive it:
\begin{Proposition}\label{prop:medial}
  The medial rule is derivable for $\{\rulecdown,\rulewdown\}$. Dually,
  the medial rule is derivable for $\{\rulecup,\rulewup\}$.
\end{Proposition}
\begin{proof}
  The medial rule is derivable as follows (or dually):
$$
\dernote {\rulecdown }{\quad.}{S\aprs{\pars{R,T},\pars{U,V}}} {
  \rootnote{\rulewdown
  }{}{S\pars{\aprs{\pars{R,T},\pars{U,V}},\aprs{\pars{R,T},\pars{U,V}}}}
  { \root{\rulewdown}{
      S\pars{\aprs{R,\pars{U,V}},\aprs{\pars{R,T},\pars{U,V}}} }{
      \root{\rulewdown}{
        S\pars{\aprs{R,U},\aprs{\pars{R,T},\pars{U,V}}} }{
        \root{\rulewdown}{ S\pars{\aprs{R,U},\aprs{T,\pars{U,V}}} }{
          \leaf{ S\pars{\aprs{R,U},\aprs{T,V}}} }}}}}
$$
\end{proof}

The medial rule has also been considered by Do\v{s}en and Petri\'c as
a composite arrow in the free bicartesian category, \cf the end of
Section~4 in \cite{DosPet02}. It is composed of four projections and a
pairing of identities (or dually) in the same way as medial is
derived using four weakenings and a contraction in the proof above.

Once we admit medial, then not only identity, cut and weakening, but
also contraction is reducible to atomic form:

\begin{Theorem}\label{thm:general}
  The rules $\ruleidown$, $\rulewdown$ and $\rulecdown$ are derivable
  for $\{\ruleaidown,\swir\}$, $\{\ruleawdown\}$ and $\{\ruleacdown,
  \rulem \}$, respectively. Dually, the rules $\ruleiup$, $\rulewup$
  and $\rulecup$ are derivable for $\{\ruleaiup,\swir\}$,
  $\{\ruleawup\}$ and $\{\ruleacup, \rulem \}$, respectively.

\end{Theorem}

\begin{proof}
  
  I will show derivability of the rules $\{\ruleidown, \rulewdown,
  \rulecdown\}$ for the respective systems. The proof of derivability
  of their co-rules is dual.
  
  Given an instance of one of the following rules:
$$
  \vcinf{\intr}{S{\pars{R,\neg R}}}{S\cons{\true}}\quad,\qquad
  \vcinf{\rulewdown}{S{\cons{R}}}{S\cons{\false}}\quad,\qquad
  \vcinf{\rulecdown}{S{\cons{R}}}{S\pars{R,R}} \quad,\qquad
  $$
construct a new derivation by structural induction on $R$:
\begin{enumerate}
\item $R$ is an atom. Then the instance of the general
rule is also an instance of its atomic form.
\item $R=\pars{P,Q}$, where $P\ne\false\ne Q$. Apply the
  induction hypothesis respectively on

$$\downsmash{
\dernote{\swir}{\quad,}{S\pars{P,Q,\aprs{        \neg P ,        \neg Q }}} {
\root   {\swir}  {S\pars{  Q,\aprs{\pars{P,\neg P},        \neg Q }}} {
\root   {\intr}  {S    \aprs{\pars{P,\neg P},\pars{Q,\neg Q}} }{
\root   {\intr}  {S        \pars{Q,\neg Q} }{
\leaf            {S\cons{\true}}                                   }}}}
}\qquad
\downsmash{
\dernote{\rulewdown}{\quad,}{S\pars{P, Q}} {
\rootnote   {\rulewdown}{}  {S \pars{\false,Q} }{
\rootnote   {=}{}  {S \pars{\false,\false} }{
\leaf            {S\cons{\false}}}}                                   }
}\qquad
\downsmash{
\dernote {\rulecdown        }{\quad.}{S\pars{P,Q}}  {
\rootnote{\rulecdown    }{}{S\pars{P,P,Q}} {
\leaf                             {S\pars{P,P,Q,Q}}}}
}
$$
\item $R=\aprs{P,Q}$, where $P\ne\true\ne Q$. Apply the
induction hypothesis respectively on 
$$\downsmash{
\dernote{\swir}{\quad,}{S\pars{\aprs{        P ,        Q },\neg P,\neg Q}} {
\root   {\swir}  {S\pars{\aprs{\pars{P,\neg P},        Q }, \neg Q}} {
\root   {\intr}  {S    \aprs{\pars{P,\neg P},\pars{Q,\neg Q}} }{
\root   {\intr}  {S        \pars{Q,\neg Q} }{
\leaf            {S\cons{\true}}                                   }}}}
}\qquad
\downsmash{
\dernote{\rulewdown}{\quad,}{S\aprs{P, Q}} {
\root   {\rulewdown} {S \aprs{\false,Q} }{
\root   {=} {S \aprs{\false,\false} }{
\leaf            {S\cons{\false}}}} }
}\qquad
\downsmash{
\dernote {\rulecdown        }{\quad.}{S\aprs{P,Q}}  {
\rootnote{\rulecdown    }{}{S\aprs{\pars{P,P},Q}} {
\root{\rulem}{S\aprs{\pars{P,P},\pars{Q,Q}}}{
\leaf{S\pars{\aprs{P,Q},\aprs{P,Q}}}}}}
}$$
\end{enumerate}

\end{proof}

We now obtain the local system $\SKS$ from $\SKSg$ by restricting
identity, cut, weakening and contraction to atomic form and adding
medial. It is shown in Figure~\ref{fig:sks}. The names of the rules
are as in system $\SKSg$, except that the atomic rules carry the
attribute \emph{atomic}, as for example in the name \emph{atomic cut}
for the rule $\ruleaiup$.

\renewcommand{\rulebox}[1]{\parbox{8em}{$$#1$$}}
\begin{figure}[htb]
  \begin{center}
    \fbox{
      \parbox{19em}{
        \rulebox{\vcinf{\ruleaidown}
          {S\pars{a,\neg a}}
          {S\cons{\true}}}
        \hfill
        \rulebox{\vcinf{\ruleaiup}
          {S\cons{\false}}
          {S\aprs{a,\neg a}}}
        $$
        \vcinf{\swir}
        {S\pars{\aprs{R,T},U }}
        {S\aprs{\pars{R,U},T}}
        $$
        $$
        \vcinf{\rulem}
        {S\aprs{\pars{R,T},\pars{U,V} }}
        {S\pars{\aprs{R,U},\aprs{T,V} }}
        $$
        \rulebox{\vcinf{\ruleawdown}
          {S\cons{a}}
          {S\cons{\false}}}
        \hfill
        \rulebox{\vcinf{\ruleawup}
          {S\cons{\true}}
          {S\cons{a}}}
        \rulebox{\vcinf{\ruleacdown}
          {S\cons{ a}}
          {S\pars{ a,a }}}
        \hfill
        \rulebox{\vcinf{\ruleacup}
          {S\aprs{ a,a }}
          {S\cons{ a}}}
        }
      }    
    \caption{System $\SKS$}
    \label{fig:sks}
  \end{center}
\end{figure}

\begin{Theorem}\label{thm:symequiv}
  System $\SKS$ and system $\SKSg$ are strongly equivalent.
\end{Theorem}

\begin{proof}
  Derivations in $\SKSg$ are translated to derivations in $\SKS$ by
  Theorem \ref{thm:general}, and vice versa by Proposition
  \ref{prop:medial}.
\end{proof}

Thus, all results obtained for the non-local system, in particular the
correspondence with the sequent calculus and admissibility of the
up-rules, also hold for the local system.  By removing the up-rules
from system $\SKS$ we obtain system $\KS$, shown in
Figure~\ref{fig:KS}.

\renewcommand{\rulebox}[1]{\mbox{$#1$}}
\begin{figure}[htb]
  \begin{center}
    \fbox{
      \parbox{30em}{\medskip
$$
        \rulebox{\vcinf{\ruleaidown}
          {S\pars{a,\neg a}}
          {S\cons{\true}}}
        \qquad\qquad
        \rulebox{\vcinf{\ruleawdown}
          {S\cons{a}}
          {S\cons{\false}}}
        \qquad\qquad
        \rulebox{\vcinf{\ruleacdown}
          {S\cons{ a}}
          {S\pars{ a,a }}}
$$
\medskip
$$
        \rulebox{\vcinf{\swir}
        {S\pars{\aprs{R,U},T }}
        {S\aprs{\pars{R,T},U}}}
      \qquad \qquad
        \rulebox{\vcinf{\rulem}
        {S\aprs{\pars{R,U},\pars{T,V} }}
        {S\pars{\aprs{R,T},\aprs{U,V} }}}
$$
        }
      }    
    \caption{System $\KS$}
    \label{fig:KS}
  \end{center}
\end{figure}

\begin{Theorem}\label{thm:equiv}
  System $\KS$ and system $\KSg$ are strongly equivalent.
\end{Theorem}

\begin{proof}
  As the proof of Theorem \ref{thm:symequiv}.
\end{proof}

In the following, I will concentrate on the local system. The
non-local rules, general identity, weakening, contraction and their
duals $\{\ruleidown, \ruleiup, \rulewdown, \rulewup, \rulecdown,
\rulecup \}$ do not belong to $\SKS$. However, I will freely use them
to denote a corresponding derivation in $\SKS$ according to Theorem
\ref{thm:general}. For example, I will use

$$
\infnote{\rulecdown}{\aprs{a,b}}{\pars{\aprs{a,b},\aprs{a,b}}}{}
$$
to denote either
$$\downsmash{ 
\dernote{\ruleacdown}{\qquad \mbox{or} \qquad}{\aprs{a,b}} {
    \root{\ruleacdown} { \aprs{\pars{a,a},b}} { { \root {\rulem}
        {\aprs{\pars{a,a},\pars{b,b}}} { \leaf
          {\pars{\aprs{a,b},\aprs{a,b}}} }}}} 
\dernote{\ruleacdown}{\qquad .}{\aprs{a,b}} {
    \root{\ruleacdown} { \aprs{a,\pars{b,b}}} { { \root {\rulem}
        {\aprs{\pars{a,a},\pars{b,b}}} { \leaf
          {\pars{\aprs{a,b},\aprs{a,b}}} }}}} 
}
$$

\subsection{Locality Through Atomicity}
In system $\SKSg$ and also in sequent systems, there is no bound on
the size of formulae that can appear as an active formula in an
instance of the contraction rule.  Implementing those systems for
proof search thus requires either duplicating formulae of unbounded
size or putting in place some complex mechanism, \eg of sharing and
copying on demand. In system $\SKS$, no rule requires duplicating
structures of unbounded size.  In fact, because no rule needs to
inspect structures of unbounded size, I call this system \emph{local}.
The atomic rules only need to duplicate, erase or compare atoms. The
switch rule involves structures of unbounded size, namely $R$, $T$ and
$U$.  But it does not require inspecting them.  To see this, consider
structures represented as binary trees in the obvious way.  Then the
switch rule can be implemented by changing the marking of two nodes
and exchanging two pointers:
$$
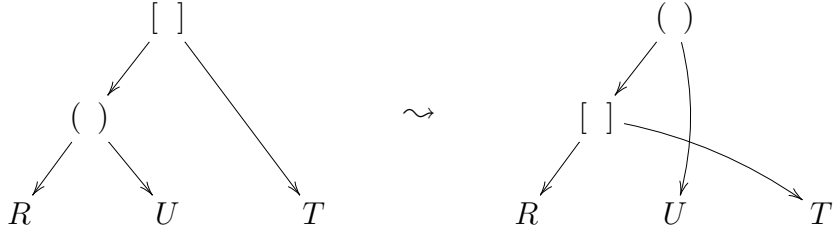
\vcenter{
\xymatrix@C=3mm@R=7mm{
      &                                & \pars{\enspace} \ar[dl] \ar[ddrrr]\\
      &   \aprs{\enspace}    \ar[dl] \ar[dr]\\ 
   R &                                &      U & & & T\\
}}
\qquad\leadsto\qquad
\vcenter{
\xymatrix@C=3mm@R=7mm{
      &                                & \aprs{\enspace} \ar[dl]
\ar@/^/[dd]\\
      &   \pars{\enspace}    \ar[dl] \ar@/^/[drrrr]\\ 
   R &                                &      U & & & T\\
}}\quad.
$$

The same technique works for medial. The equations are local as well,
including the De Morgan laws.  However, since the rules in $\SKS$
introduce negation only on atoms, it is even possible to restrict
negation to atoms from the beginning, as is customary in the one-sided
sequent calculus, and drop the equations for negation entirely.

The concept of locality depends on the representation of structures.
Rules that are local for one representation may not be local when
another representation is used.  For example, the switch rule is local
when structures are represented as trees, but it is not local when
structures are represented as strings.

One motivation for locality is to simplify distributed implementation
of an inference system. Of course, locality by itself still makes no
distributed implementation. There are tasks to accomplish in an
implementation of an inference system that in general require a global
view on structures, for example matching a rule, \ie finding a redex.
There should also be some mechanism for backtracking. I do not see how
these problems can be approached within a proof-theoretic system with
properties like cut elimination.  However, the application of a rule,
\ie producing the contractum from the redex, is achieved locally in
system $\SKS$. For that reason I believe that it lends itself more
easily to distributed implementation than other systems.

\subsection{Decomposition of Derivations}\label{sec:decomp}

Derivations can be arranged into consecutive phases such that each
phase uses only certain rules. We call this property
\emph{decomposition}.  Decomposition theorems thus provide normal
forms for derivations.  A classic example of a decomposition theorem
in the sequent calculus is proving Herbrand's Theorem by decomposing
a proof tree into a bottom phase with contraction and quantifier rules
and a top phase with propositional rules only. The three decomposition
theorems presented here state the possibility of pushing all instances
of a certain rule to the top and all instances of its dual rule to the
bottom of a derivation. Except for the first, these decomposition
theorems do not have analogues in the sequent calculus.

The proofs in this section use \emph{permutation} of rules.

\begin{Definition}
      A rule $\rho$ \emph{permutes over} a rule $\pi$ (or $\pi$
    \emph{permutes under} $\rho$) if for every derivation
    $\vcenter{\dernote{\rho}{}{R}{\root{\pi}{U}{\leaf{T}}}}$ there is
    a derivation
    $\vcenter{\dernote{\pi}{}{R}{\root{\rho}{V}{\leaf{T}}}}$ for some
    structure $V$.  
\end{Definition}

\begin{Lemma}\label{lem:permac}
  The rule $\ruleacdown$ permutes under the rules $\ruleawdown$,
  $\ruleaidown$, $\rules$ and $\rulem$. Dually, the rule $\ruleacup$
  permutes over the rules $\ruleawup$, $\ruleaiup$, $\rules$ and
  $\rulem$.
\end{Lemma}

\begin{proof}
  Given an instance of $\ruleacdown$ above an instance of a rule
  $\rho\in \{\ruleawdown, \ruleaidown, \rules, \rulem\}$, the redex of
  $\ruleacdown$ can be a substructure of the context of $\rho$. Then
  we permute as follows:
  $$
  \vcenter{\dernote{\rho}{}{S\cons{R}}{\root{\ruleacdown}{S\cons{U}}{\leaf{S'\cons{U}}}}}
  \qquad \leadsto \qquad
  \vcenter{\dernote{\ruleacdown}{}{S\cons{R}}{\root{\rho}{S'\cons{R}}{\leaf{S'\cons{U}}}}}
\quad .
$$
The only other possibility occurs in case that $\rho$ is $\rules$
or $\rulem$: the redex of $\ruleacdown$ can be a substructure of the
contractum of $\rho$. Then we permute as in the following example of a
switch rule, where $T\cons{\enspace}$ is a structure context:
$$
\vcenter{
  \dernote{\rules}{}{S\pars{R,\aprs{T\cons{a},U}}}{\root{\ruleacdown}{S\aprs{\pars{R,T\cons{a}},U}}{\leaf{S\aprs{\pars{R,T\pars{a,a}},U}
      }}}} \qquad \leadsto \qquad \vcenter{
  \dernote{\ruleacdown}{}{S\pars{R,\aprs{T\cons{a},U}}}{\root{\rules}{S\aprs{\pars{R,T\pars{a,a}},U}}{\leaf{S\aprs{\pars{R,T\pars{a,a}},U}
      }}} } \quad .
$$
(And dually for $\ruleacup$.)
\end{proof}

\begin{Lemma}\label{lem:permaw}
  The rule $\ruleawdown$ permutes under the rules $\ruleaidown$,
  $\rules$ and $\rulem$. Dually, the rule $\ruleawup$ permutes over
  the rules $\ruleaiup$, $\rules$ and $\rulem$.
\end{Lemma}
\begin{proof}
  Similar to the proof of Lemma \ref{lem:permac}.
\end{proof}

We now turn to the decomposition results.

\subsubsection*{Separating Identity and Cut}

Given that in system $\SKS$ identity is a rule, not an axiom as in the
sequent calculus, a natural question to ask is whether the
applications of the identity rule can be restricted to the top of a
derivation. For proofs, this question is already answered positively
by Theorem \ref{thm:semcutelim}. It turns out that it is also true for
derivations in general. Because of the duality between $\ruleaidown$
and $\ruleaiup$ we can also push the cuts to the bottom of a
derivation.  While this can be obtained in the sequent calculus (using
cut elimination), it can not be done with a simple permutation
argument as we do.

The following rules are called \emph{super switch down} and \emph{super
  switch up}:
$$
\vcinf{\rulessdown} {S\pars{R,T\cons{\false}}} {S\cons{T\cons{R}}}
\qquad{\rm and}\qquad \vcinf{\rulessup} {S\cons{T\cons{R}}}
{S\aprs{R,T\cons{\true}}} \quad .
  $$

\begin{Lemma}\label{lem:superswitch}
  The rule $\rulessdown$ is derivable for $\{\rules\}$. Dually, the
  rule $\rulessup$ is derivable for $\{\rules\}$.
\end{Lemma}

\begin{proof}
  We prove this for $\rulessup$ by structural induction on
  $T\cons{\enspace}$. The proof for $\rulessdown$ is dual.
\begin{enumerate}
\item $T\cons{\enspace}$ is empty. Then premise and conclusion of the
  given instance of $\rulessup$ coincide, the rule instance can
  be removed.
  
\item $T\cons{\enspace}=\pars{U,V\cons{\enspace}}$, where
  $U\ne\false$. Apply the induction hypothesis on
$$\downsmash{ \dernote {\rulessup }{}{S\pars{U,V\cons{R}}} {
     \root{\rules}{S\pars{U,\aprs{R,V\cons{\true}}}}{
         \leaf{S\aprs{R,\pars{U,V\cons{\true}}}}}} }
\quad .$$

\item $T\cons{\enspace}=\aprs{U,V\cons{\enspace}}$, where
  $U\ne\true$. Apply the induction hypothesis on
$$
\downsmash{ 
\infnote{\rulessup}{ S\aprs{U,V\cons{R}}}{S\aprs{U,R,V\cons{\true} }}
{\quad , } }
$$
the redex being $V\cons{R}$.
\end{enumerate}
\end{proof}

The following rules are called \emph{shallow} atomic identity and cut:
$$
\vcinf{\ruleaisdown}
        {\aprs{S,\pars{a,\neg a}}}
        {S}
\qquad{\rm and}\qquad
\vcinf{\ruleaisup}
        {S}
        {\pars{S,\aprs{a,\neg a}}}
\quad .
$$

\begin{Lemma}\label{lem:shallowint}
  The rule $\ruleaidown$ is derivable for
  $\{\ruleaisdown,\swir\}$. Dually, the rule $\ruleaiup$ is derivable
  for $\{\ruleaisup,\swir\}$.
\end{Lemma}
\begin{proof}
  An instance of $\ruleaidown$ can be replaced by an instance of
  $\ruleaisdown$ followed by an instance of $\rulessup$, which is
  derivable for $\{\rules\}$. (And dually for $\ruleaiup$.)
\end{proof}

\begin{Theorem}\label{thm:sepint}
\theoremnl
    For every derivation $\vcenter{\xy\xygraph{[]!{0;<1pc,0pc>:} {T}
        -@{=}^<>(.5){\strut\SKS}[dd] {R}
        }\endxy}$  there is a derivation
    $\vcenter{\xy\xygraph{[]!{0;<1pc,0pc>:} {T}
        -@{=}^<>(.5){\strut\{\ruleaidown\}}[dd] {V}
        -@{=}^<>(.5){\strut\SKS \, \setminus \{\ruleaidown, \ruleaiup\}}[dd] {U}
        -@{=}^<>(.5){\strut\{\ruleaiup\}} [dd] {R}
}\endxy}$.
  
\end{Theorem}

\begin{proof}
  By Lemma \ref{lem:shallowint} we can reduce atomic identities to
  shallow atomic identities and the same for the cuts. It is easy to
  check that the rule $\ruleaisdown$ permutes over every rule in
  $\SKS$ and the rule $\ruleaisup$ permutes under every rule in
  $\SKS$. Instances of $\ruleaisdown$ and $\ruleaisup$ are instances
  of $\ruleaidown$ and $\ruleaiup$, respectively.  
\end{proof}

\subsubsection*{Separating Contraction}

Contraction allows the repeated use of a structure in a proof by
allowing us to copy it at will. It should be possible to copy
everything needed in the beginning, and then go on with the proof
without ever having to copy again. This intuition is made precise by
the following theorem and holds for system $\SKS$. There is no such
result for the sequent calculus \cite{BruRC02}. There are sequent
systems for classical propositional logic that do not have an explicit
contraction rule, however, since they involve context sharing,
contraction is built into the logical rules and is used throughout the
proof.

\begin{Theorem}\label{thm:sepac1}\theoremnl
For every proof $\vcenter{\xy\xygraph{[]!{0;<1pc,0pc>:}
      {}*=<0pt>{}:@{|=}^<>(.5){\strut\KS}[dd] {S} }\endxy}$ there is
  a proof $\vcenter{\xy\xygraph{[]!{0;<1pc,0pc>:}
      {}*=<0pt>{}:@{|=}^<>(.5){\strut\KS\setminus\{\ruleacdown\}}[dd]
      {S'} -@{=}^<>(.5){\strut\{\ruleacdown\}} [dd] {S} }\endxy}$.
\end{Theorem}

\begin{proof}
  Using Lemma \ref{lem:permac}, permute down all instances of
  $\ruleacdown$, starting with the bottommost.
\end{proof}

This result is extended to the symmetric system as follows: 

\begin{Theorem}\label{thm:sepac2}\theoremnl
  For every derivation $\vcenter{\xy\xygraph{[]!{0;<1pc,0pc>:} {T}
      -@{=}^<>(.5){\strut\SKS}[dd] {R} }\endxy}$ there is a derivation
  $\vcenter{\xy\xygraph{[]!{0;<1pc,0pc>:} {T}
      -@{=}^<>(.5){\strut\{\ruleacup\}}[dd] {V}
      -@{=}^<>(.5){\strut\SKS \, \setminus \{\ruleacdown,
        \ruleacup\}}[dd] {U} -@{=}^<>(.5){\strut\{\ruleacdown\}} [dd]
      {R} }\endxy}$.
\end{Theorem}

\begin{proof}
  Consider the following derivations that can be obtained:
  
  $$
  \vcenter{\xy\xygraph{[]!{0;<1pc,0pc>:} {T}
      -@{=}^<>(.5){\strut\SKS}_<>(.5){\strut}[ddd] {R} }\endxy}
  \stackrel{1}{\leadsto} \vcenter{\xy\xygraph{[]!{0;<1pc,0pc>:}
      {\inf{\ruleidown}{\pars{\neg T, T}}{\true}\hspace{2.7ex}}
      -@{=}^<>(.5){\strut\SKS} _<>(.5){\strut}[ddddd] {\pars{ \neg T,
          R}} }\endxy} \stackrel{2}{\leadsto}
  \vcenter{\xy\xygraph{[]!{0;<1pc,0pc>:} {\true}
      -@{=}^<>(.5){\strut\KS} [ddd] {\pars{\neg T,R}} }\endxy}
  \stackrel{3}{\leadsto} \vcenter{\xy\xygraph{[]!{0;<1pc,0pc>:}
      {\aprs{T,\true}} -@{=}^<>(.5){\strut\KS} [ddddd]
      {\dernote{\aintr}{}{R}{ \root{\swir}{\pars{R,\aprs{T,\neg T}}}{
            \leaf{\aprs{T,\pars{\neg T,R}} }}}} }\endxy}
  \stackrel{4}{\leadsto} \vcenter{\xy\xygraph{[]!{0;<1pc,0pc>:}
      {\aprs{T,\true}} -@{=}^<>(.5){\strut\KS\setminus\{\ruleacdown\}}
      [ddd]{\pars{R',\aprs{T,\neg T'} }}
      -@{=}^<>(.5){\strut\{\ruleacdown\}} [dddd]
      {\infnote{\aintr}{R'}{\pars{R',\aprs{T,\neg T}}}{}}
      -@{=}^<>(.5){\strut\{\ruleacdown\}} [dddd] {R} } \endxy}
\quad .
$$
\begin{enumerate}
\item Put the derivation into the context $\pars{\neg
    T,\cons{\enspace}}$. On top of the resulting derivation, apply an
  $\ruleidown$ to obtain a proof.
\item Apply cut elimination.
\item Put the proof into the context $\aprs{T,\cons{\enspace}}$. At
  the bottom of the resulting derivation, apply a switch and a cut to
  obtain a derivation from $T$ to $R$.
\item All instances of $\ruleacdown$ are permuted down as far as
  possible. There are two kinds of instances: those that duplicate
  atoms coming from $R$ and those that duplicate atoms coming from
  $\neg T$.  The first kind, starting with the bottom-most instance,
  can be permuted down all the way to the bottom of the derivation.
  The second kind, also starting with the bottom-most instance, can be
  permuted down until they meet the cut.
\end{enumerate}

Now, starting with the bottom-most $\ruleacdown$ that is above the cut,
we apply the transformation
$$
\dernote{\aintr}{\quad \leadsto \quad}{S\cons{\false}}{
            \root{\ruleacdown}{S\aprs{U\cons{a},\neg U\cons{\neg a}}}{
                \leaf{S\aprs{U\cons{a},\neg U\pars{\neg a, \neg a}}}}}
\dernote{\aintr}{}{S\cons{\false}}{
            \root{\ruleacup}{S\aprs{U\aprs{a,a},\neg U\pars{\neg a,
\neg a}}}{
                \leaf{S\aprs{U\cons{a},\neg U\pars{\neg a, \neg a}}}}}
            $$
            and permute up the resulting instance of $\ruleacup$
            all the way to the top of the derivation. This is possible
            because no rule in the derivation above changes $T$.
            Proceed inductively with the remaining instances of
            $\ruleacdown$ above the cut. The resulting derivation has
            the desired shape.   
\end{proof}

\subsubsection*{Separating Weakening}

In the sequent calculus, one usually can push up to the top of the
proof all the instances of weakening, or to the same effect, build
weakening into the identity axiom:
$$
\inf{}{A, \Phi \vdash  A, \Psi}{} \quad .
$$
The same \emph{lazy} way of applying weakening can be done in
system $\SKS$, \cf Theorem \ref{thm:semcutelim}. However, while a
proof in which all weakenings occur at the top is certainly more
``normal'' than a proof in which weakenings are scattered all over,
this is hardly an interesting normal form. In system $\SKS$ something
more interesting can be done: applying weakening in an \emph{eager}
way.

\begin{Theorem}\label{thm:sepaw1}\theoremnl
  For every proof $\vcenter{\xy\xygraph{[]!{0;<1pc,0pc>:}
      {}*=<0pt>{}:@{|=}^<>(.5){\strut\KS}[dd] {S} }\endxy}$ there is
  a proof $\vcenter{\xy\xygraph{[]!{0;<1pc,0pc>:}
      {}*=<0pt>{}:@{|=}^<>(.5){\strut\KS\setminus\{\ruleawdown\}}[dd]
      {S'} -@{=}^<>(.5){\strut\{\ruleawdown\}} [dd] {S} }\endxy}$.
\end{Theorem}

\begin{proof}
  Permute down all instances of $\ruleawdown$, starting with the
  bottommost. This is done by using Lemma \ref{lem:permaw} and the
  following transformation:
  $$
  \dernote{\ruleacdown}{\quad \leadsto
    \quad}{S\cons{a}}{
    \root{\ruleawdown}{S\pars{a,a}}{
      \leaf{S\pars{a,\false}}}}
  \dernote{=}{\quad .}{S\cons{a}}{ \leaf{S\pars{a,\false}}}
  $$
\end{proof}

Weakening loses information: when deducing $a \vee b$ from $a$, the
information that $a$ holds is lost.  Given a proof with weakenings, do
they lose information that we would like to keep? Can we obtain a
proof of a stronger statement by removing them?  The theorem above
gives an affirmative answer to that question: given a proof of $S$, it
exhibits a weakening-free proof of a structure $S'$, from which $S$
trivially follows by weakenings.

\begin{Notation}\hspace{0mm}\\[1ex]
  A derivation $\vcenter{\xy\xygraph{[]!{0;<2pc,0pc>:}
      {T}-@{=}^<>(.5){\strut\{\rho\}} _<>(.5){\strut}[d] {R} }\endxy}
  $ of $n$ instances of the rule $\rho$ is denoted by
  $\vcenter{\infnote{\rho^n}{\; R \; }{\; T \; }{\; .}}  $
\end{Notation}

This result extends to the symmetric case, \ie to derivations in
$\SKS$:

\begin{Theorem}\label{thm:sepaw2}\theoremnl
For every derivation $\vcenter{\xy\xygraph{[]!{0;<1pc,0pc>:} {T}
    -@{=}^<>(.5){\strut\SKS}[dd] {R} }\endxy}$ there is a derivation
$\vcenter{\xy\xygraph{[]!{0;<1pc,0pc>:} {T}
    -@{=}^<>(.5){\strut\{\ruleawup\}}[dd] {V} -@{=}^<>(.5){\strut\SKS
      \, \setminus \{\ruleawdown, \ruleawup\}}[dd] {U}
    -@{=}^<>(.5){\strut\{\ruleawdown\}} [dd] {R} }\endxy}$.
\end{Theorem}

\begin{proof}
  
  There is an algorithm that produces a derivation of the desired
  shape. It consists of two steps: 1) pushing up all instances of
  $\ruleawup$ and 2) pushing down all instances of $\ruleawdown$.
  Those two steps are repeated until the derivation has the desired
  shape. An instance of $\ruleawup$ that is pushed up can turn into an
  instance of $\ruleawdown$ when meeting an instance of $\ruleaidown$,
  and the dual case can also happen. However, the process terminates,
  since each step that does not produce the desired shape strictly
  decreases the combined number of instances of $\ruleaidown$ and
  $\ruleaiup$.
  
  In the following, the process of pushing up instances of $\ruleawup$
  is shown, the process of pushing down instances of $\ruleawdown$ is
  dual. 
  
  An instance of $\ruleawup$ is a special case of a derivation
  consisting of $n$ instances of $\ruleawup$ and is moved up as
  such, starting with the topmost.  In addition to the cases treated
  in Lemma \ref{lem:permaw} there are the following cases:

  $$
  \dernote{\ruleawup^n}{\quad \leadsto
    \quad}{S\aprs{\true,a}}{
    \root{\ruleacup}{S'\aprs{a,a}}{
      \leaf{S'\cons{a}}}}
  \dernote{=}{}{S\aprs{\true,a}}{
    \root{\ruleawup^{n-1}}{S\cons{a}}{
      \leaf{S'\cons{a}}}}
  $$

$$
  \dernote{\ruleawup^n}{\quad \leadsto
    \quad}{S\cons{\true}}{
    \root{\ruleacdown}{S'\cons{a}}{
      \leaf{S'\pars{a,a}}}}
  \dernote{=}{}{S\cons{\true}}{
    \root{\ruleawup^{n+1}}{S\pars{\true,\true}}{
      \leaf{S'\pars{a,a}}}}
  $$

  $$\dernote{\ruleawup^n}{\quad \leadsto \quad}{S\cons{\true}}{
    \root{\ruleawdown}{S'\cons{a}}{ \leaf{S'\cons{\false}}}} \dernote{
    = }{}{ S\cons{\true} }{
    \root{\rules}{S\pars{\true,\aprs{\true,\false}}}{
      \root{=}{S\aprs{\pars{\true,\true},\false} }{
        \root{\ruleawup^{n-1}}{ S\cons{\false} }{ \leaf{
            S'\cons{\false} }}}}}
$$
  
$$\dernote{ \ruleawup^n }{\quad \leadsto \quad}{ S\pars{\true, \neg
    a}}{ \root{ \ruleaidown }{ S'\pars{a, \neg a}}{ \leaf{
      S'\cons{\true}}}} \dernote{ \ruleawdown }{}{ S\pars{\true, \neg
    a}}{ \root{=}{ S\pars{\true, \false}}{
    \root{\ruleawup^{n-1}}{S\cons{\true}}{ \leaf{ S'\cons{\true}}}} }
  $$
\end{proof}

\subsubsection*{Separating all Atomic Rules}

Decomposition results can be applied consecutively. Here, all rules
that deal with atoms, namely $\ruleaidown$, $\ruleacdown$,
$\ruleawdown$ and their duals, are separated from the rules that deal
with the connectives, namely $\rules$ and $\rulem$:

\begin{Theorem}\label{thm:sepall}\theoremnl
  For every derivation $\vcenter{\xy\xygraph{[]!{0;<1pc,0pc>:} {T}
      -@{=}^<>(.5){\strut\SKS}[dd] {R} }\endxy}$ there is a derivation
  $\vcenter{\xy\xygraph{[]!{0;<1pc,0pc>:} {T}
      -@{=}^<>(.5){\strut\{\ruleacup\}}[dd] {T_1}
      -@{=}^<>(.5){\strut\{\ruleawup\}}[dd] {T_2}
      -@{=}^<>(.5){\strut\{\ruleaidown\}}[dd] {T_3}
      -@{=}^<>(.5){\strut\{\rules,\rulem\}}[dd] {R_3}
      -@{=}^<>(.5){\strut\{\ruleaiup\}}[dd] {R_2}
      -@{=}^<>(.5){\strut\{\ruleawdown\}}[dd] {R_1}
      -@{=}^<>(.5){\strut\{\ruleacdown\}} [dd] {R} }\endxy} \qquad.$
\end{Theorem}
\begin{proof}
  We first push contractions to the outside, using Theorem
  \ref{thm:sepac2}. In the contraction-free part of the obtained
  derivation, we push weakening to the outside, using the procedure
  from the proof of Theorem \ref{thm:sepaw2}, which does not introduce
  new instances of contraction. In the contraction- and weakening-free
  part of the resulting derivation we then separate out identity and
  cut by applying the procedure from the proof of Theorem
  \ref{thm:sepint}, which neither introduces new contractions nor
  weakenings.
\end{proof}

\section{Predicate Logic}\label{sec:predicate}

\begin{Definition}
  \emph{Variables} are denoted by $x$ and $y$.  Terms are defined as
  usual in first-order predicate logic. \emph{Atoms}, denoted by $a,
  b$, etc., are expressions of the form $p(t_1,\dots ,t_n)$, where $p$
  is a \emph{predicate symbol} of \emph{arity} $n$ and $t_1,\dots
  ,t_n$ are terms.  The negation of an atom is again an atom.  The
  \emph{structures} of the language $\KSpred$ are generated by the
  following grammar, which is the one for the propositional case
  extended by existential and universal quantifier:
  $$
  S \grammareq \false \mid \true \mid a \mid
  \pars{\,\underbrace{S,\dots,S}_{{}>0}\,} \mid
  \aprs{\,\underbrace{S,\dots,S}_{{}>0}\,} \mid \neg S \mid \exists x
  S \mid \forall x S\quad .
  $$
\end{Definition}

\begin{Definition}
  Structures are \emph{equivalent} modulo the smallest equivalence
  relation induced by the equations in Fig.~\ref{fig:Equations}
  extended by the following equations:
        $$
         \begin{array}[h]{lcl}
        \mbox{\smalltitle{Variable Renaming}} \quad &
           \begin{array}[h]{c}
             \forall x R = \forall y R\subst{x/y} \\
           \exists x R = \exists y R\subst{x/y} \\
           \end{array}
        & \mbox{if $y$ is not free in $R$}\\[3ex]
        \mbox{\smalltitle{Vacuous Quantifier}} &
          \forall y R = \exists y R = R     & \mbox{if $y$ is not free in $R$}\\[1ex]
        \mbox{\smalltitle{Negation}} &
        \begin{array}[h]{l}
          \overline{\exists x R} = \forall x \neg R\\
          \overline{\forall x R} = \exists x \neg R\\
        \end{array}
        \end{array}
        $$
\end{Definition}

\begin{Definition}
  The notions of \emph{structure context} and \emph{substructure} are
  defined in the same way as in the propositional case. A structure of
  language $\KSpred$ is in \emph{normal form} if negation occurs only
  on atoms, and extra units, connectives and quantifiers are removed
  using the laws for units, associativity and vacuous quantifier.
\end{Definition}

As in the
propositional case, we in general consider structures to be in normal
form and do not distinguish between equivalent structures.

\renewcommand{\rulebox}[1]{\parbox{8em}{$$#1$$}}
\begin{figure}[t]
  \begin{center}
    \fbox{
      \parbox{.8\textwidth}{
        \hfill
        \rulebox{\vcinf{\ruleidown}
          {S\pars{R,\neg R}}
          {S\cons{\true}}}
        \qquad
        \rulebox{\vcinf{\ruleiup}
          {S\cons{\false}}
          {S\aprs{R,\neg R}}}
        \hfill
        \hspace{0mm}
\\ %\vspace{0mm} \dotfill \vspace{0mm} \\
       \vspace{0mm}
       \hfill
       \rulebox{\vcinf{\swir}
        {S\pars{\aprs{R,T},U}}
        {S\aprs{\pars{R,U},T}}}
        \hfill
        \vspace{-10mm}
\\
       \rulebox{\vcinf{\ruleudown}
          {S\pars{\forall x R, \exists x T}}
          {S\cons{\forall x \pars{R,T}}}}
        \hfill
        \rulebox{\vcinf{\ruleuup}
          {S\cons{\exists x \aprs{R,T}}}
          {S\aprs{\exists x R, \forall x T}}} \vspace{5mm}
\\
        \vspace{0mm}
        \hfill
        \rulebox{\vcinf{\rulewdown}
          {S\cons{R}}
          {S\cons{\false}}}
        \qquad
        \rulebox{\vcinf{\rulewup}
          {S\cons{\true}}
          {S\cons{R}}}
        \hfill 
        \vspace{0mm}
\\ 
        \vspace{0mm}
        \hfill
        \rulebox{\vcinf{\rulecdown}
          {S\cons{ R}}
          {S\pars{ R,R }}}
        \qquad
        \rulebox{\vcinf{\rulecup}
          {S\aprs{ R,R }}
          {S\cons{ R}}}
        \hfill  
        \vspace{0mm}
\\
        \rulebox{\vcinf{\rulendown}
          {S\cons{\exists x R}}
          {S\cons{R\subst{x/t}}}}
        \hfill
        \rulebox{\vcinf{\rulenup}
          {S\cons{R\subst{x/t}}}
          {S\cons{\forall x R}}}
         }
      }    
    \caption{System $\SKSgpred$}
    \label{fig:sksgq}
  \end{center}
\end{figure}

\subsection{A Symmetric System}

The rules of system $\SKSgpred$, a symmetric system for predicate
logic, are shown in Figure~\ref{fig:sksgq}. The first and last column
show the rules that deal with quantifiers, in the middle there are the
rules for the propositional fragment. The rules $\ruleudown$ and
$\ruleuup$ were obtained by Guglielmi. They follow a scheme or recipe
\cite{GugRecipe02}, which also yields the switch rule and ensures
atomicity of cut and identity not only for classical logic but also
for several other logics. The $\ruleudown$ rule corresponds to the
$\sruleRA$ rule in $\GSone$, shown in Figure \ref{fig:quantifiers}.
Because of the equational theory, we can equivalently replace it by
$$
\rulebox{\vcinf{\ruleudown} {S\pars{\forall x R, T}}
  {S\cons{\forall x \pars{R,T}}}} \quad \mbox{if $x$ is not free in
  $T$.}
$$
In the sequent calculus, going up, the $\sruleRA$ rule removes a
universal quantifier from a formula to allow other rules to access
this formula.  In system $\SKSgpred$, inference rules apply deep
inside formulae, so there is no need to remove the quantifier: it can
be moved out of the way using the rule $\ruleudown$ and it vanishes
once the proof is complete because of the equation $\forall x
\tinyspace \true = \true$.  As a result, the premise of the
$\ruleudown$ rule implies its conclusion, which is not true for the
$\sruleRA$ rule of the sequent calculus. The $\sruleRA$ rule is the
only rule in $\GSone$ with such bad behaviour.  In all the rules that
I presented in the calculus of structures the premise implies the
conclusion.

The rule $\rulendown$ corresponds to $\sruleRE$.  As usual, the
substitution operation requires $t$ to be free for $x$ in $R$:
quantifiers in $R$ do not capture variables in $t$.  The term $t$ is
not required to be free for $x$ in $S\cons{R}$: quantifiers in $S$ may
capture variables in $t$.

\subsection{Correspondence to the Sequent Calculus}

We extend the translations between $\SKSg$ and $\GSonep$ to
translations between $\SKSgpred$ and $\GSone$. System $\GSone$ is
system $\GSonep$ extended by the rules shown in
Figure~\ref{fig:quantifiers}.

\renewcommand{\rulebox}[1]{\parbox{8em}{$$#1$$}}
\begin{figure}[tb]
\begin{center}
  \fbox{\parbox{22em}{
    \begin{center}
      \rulebox{\vcinf{\sruleRE}{\sqn{\Phi, \exists x
            A}}{\sqn{\Phi,A\subst{x/t}}}}\qquad
      \rulebox{\vcinf{\sruleRA}{\sqn{\Phi, \forall x
            A}}{\sqn{\Phi,A\subst{x/y}}}}\\
      Proviso: $y$ is not free in the conclusion of $\sruleRA$.
\end{center}
} }
\caption{Quantifier rules of $\GSone$}
  \label{fig:quantifiers}
\end{center}
\end{figure}

The functions $\ls{\enspace . \enspace}$ and $\lg{\enspace
. \enspace}$ are extended in the obvious way:
$$
\begin{array}{ccc}
\ls{\exists x A}  & = & \exists x \ls A      \\[1ex] 
\ls{\forall  x A}  & = & \forall x \ls A      \\
\end{array}
\qquad \mbox{and} \qquad
\begin{array}{ccc}
\lg{\exists x S}  & = & \exists x \lg S      \\[1ex] 
\lg{\forall  x S}  & = & \forall x \lg S      \\ 
\end{array}
\quad .
$$

\subsubsection*{From the Sequent Calculus to the Calculus of
Structures}

\newcommand{\KSpredcut}{\SKSgpred\,\setminus\,\{\rulewup,\rulecup,\ruleuup,\rulenup\}}
\begin{Theorem}
  \theoremnl
  For every derivation
  $\vcenter{\Derivation{\Sigma_1}{\enspace\cdots\enspace}{\Sigma_h}\Sigma{}}$
  in $\GSone+\sruleCut$, in which the free variables in the premises that are
  introduced by $\sruleRA$ instances are $x_1,\dots, x_n$, there exists
  a derivation $\vcenter{ \xy \xygraph{[]!{0;<3pc,0pc>:} {\forall
        x_1\dots\forall x_n
        \aprs{{\ls{\Sigma_1}},\dots,\ls{\Sigma_h}}}-@{=}^<>(.5){\strut\KSpredcut}[d]
      { \ls{\Sigma} } } \endxy }$ with the same number of cuts.
\end{Theorem}

\begin{proof}
  The proof is an extension of the proof of Theorem \ref{thm:GS1toSKS}.
  There are two more inductive cases, one for $\sruleRE$, which is
  easily translated into an $\rulendown$, and one for $\sruleRA$,
  which is shown here:

$$\vcenter{ \dernote{\sruleRA}{}{\sqn{\Phi, \forall x
A}}{
    \Derivationleaf {\Sigma_1} {\enspace\cdots\enspace} {\Sigma_{h'}}
    {\sqn{\Phi, A\subst{x/y}}} {} }  } \qquad .
$$  

By induction hypothesis we have the derivation 
$$
\vcenter{\xy \xygraph{[]!{0;<4pc,0pc>:} {\forall x_1\dots\forall
      x_{n'}\aprs{{\ls{\Sigma_1}},\dots,\ls{\Sigma_{h'}}}}-@{=}_{\Delta}^<>(.5){\strut\KSpredcut}[d]
    { \pars{\ls{\Phi},\ls{A\subst{x/y}}} } } \endxy} \qquad ,
$$
from which we build
$$
\vcenter{\xy\xygraph{[]!{0;<6.5pc,0pc>:} {\forall y\forall
      x_1\dots\forall x_{n'}\aprs{\ls{\Sigma_1},\dots,\ls{\Sigma_{h'}}}}
    -@{=}_<>(.5){\forall y\cons{\Delta}}^<>(.5){\strut\KSpredcut}[d]
    { \dernote {=}{}{\pars{\ls{\Phi}, \ls{\forall x A}}} {
        \rootnote{=}{}{ \pars{ \ls{\Phi}, \ls{\forall y A\subst{x/y}}
          } } { \root{\ruleudown}{ \pars{ \exists y \ls{\Phi},
              \ls{\forall y A\subst{x/y}} } }{\leaf { \forall
              y\pars{\ls{\Phi}, \ls{A\subst{x/y}} } }} }} } }\endxy}
\qquad , 
$$
where in the lower instance of the equivalence rule $y$ is not free in
$\ls{\forall x A}$ and in the upper instance of the equivalence rule $y$ is
not free in $\ls\Phi$: both due to the proviso of the $\sruleRA$ rule.
\end{proof}

\begin{Corollary}\label{cor:gs1tosks}
  If a sequent $\Sigma$ has a proof in $\GSone+\sruleCut$ then $\ls{\Sigma}$ has
  a proof in $\KSpredcut$.
\end{Corollary}

\begin{Corollary}\label{cor:gs1tosks2}
  If a sequent $\Sigma$ has a proof in $\GSone$ then $\ls{\Sigma}$ has
  a proof in
  $\SKSgpred\,\setminus\,\{\ruleiup,\rulewup,\rulecup,\ruleuup,\rulenup\}$
  .
\end{Corollary}

\subsubsection*{From the Calculus of Structures to the Sequent Calculus}

\begin{Lemma}\label{lem:contextpred}
  For every two formulae $A,B$ and every formula context
  $C\cons{\enspace}$ there exists a derivation
  $\vcenter{\Derivation{\quad}{\sqn{A,\neg
        B}}{\quad}{\sqn{C\cons{A},\overline{C\cons{B}}}}{}}$ in
  $\GSone$.
\end{Lemma}

\begin{proof}
  There are two cases needed in addition to the proof of Lemma
  \ref{lem:context}: $C\cons{\enspace}=\exists x C'\cons{\enspace}$
  and $C\cons{\enspace}=\forall x C'\cons{\enspace}$. The first case
  is shown here, the second is similar:
$$
\dernote{\sruleRA}{\quad ,}{\sqn{\exists x C'\cons{A},\forall x
    \overline{C'\cons{B}}}}{\root{\sruleRE}{ \sqn{\exists x
      C'\cons{A},\overline{C'\cons{B}}}
  }{\leaf{\Derivation{\quad}{\sqn{A,\neg
          B}}{\quad}{\sqn{C'\cons{A},\overline{C'\cons{B}}}}\Delta
    }}}
  $$
  where $\Delta$ exists by induction hypothesis. 
\end{proof}

\begin{Theorem}\label{thm:SKSpredtoGS1}
  For every derivation $\vcenter{ \xy \xygraph{[]!{0;<3pc,0pc>:}
      {Q}-@{=}_{}^<>(.5){\strut\SKSgpred}[d] {P} } \endxy }$ there
  exists a derivation
  $\vcenter{\Derivation{\quad\quad}{\sqn{\lg{Q}}}{\quad\quad}{\sqn{\lg{P}}}{}}$
  in $\GSone+\sruleCut$.
\end{Theorem}

\begin{proof} 
  The proof is an extension of the proof of Theorem
  \ref{thm:SKStoGS1}. The base cases are the same, in the inductive
  cases the existence of $\Delta_1$ follows from Lemma
  \ref{lem:contextpred}. Corresponding to the rules for quantifiers,
  there are four additional inductive cases.  For the rules
$$
       \vcinf{\ruleudown}
          {S\pars{\forall x R, \exists x T}}
          {S\cons{\forall x \pars{R,T}}}
\qquad \mbox{and} \qquad
       \vcinf{\ruleuup}
          {S\cons{\exists x \aprs{R,T}}}
          {S\aprs{\exists x R, \forall x T}}
$$
we have the proofs
$$
\vcenter{ \dernote{\sruleRv}{}{\sqn{\forall x R \vee \exists x T,
      \exists x(\neg R \wedge \neg T)}}{ \root{\sruleRA}{\sqn{\forall
        x R, \exists x T, \exists x(\neg R \wedge \neg T)}}{
      \root{\sruleRE}{\sqn{R, \exists x T, \exists x(\neg R \wedge
          \neg T)}}{ \root{\sruleRE}{\sqn{ R, \exists x T, \neg R
            \wedge \neg T}}{ \rroot{\sruleRand}{\sqn{ R, T, \neg R
              \wedge \neg T}}{ \root{\sruleAx}{\sqn{R,\neg R}}{
              \leaf{}}}{ \root{\sruleAx}{\sqn{T, \neg T}}{
              \leaf{}}}}}}} } \qquad \mbox{and} \qquad \vcenter{
  \dernote{\sruleRv}{}{\sqn{\exists x(R \wedge T), \forall x
      \neg R \vee \exists x \neg T}}{ \root{\sruleRA}{\sqn{\exists x(R
        \wedge T), \forall x \neg R, \exists x \neg T}}{
      \root{\sruleRE}{\sqn{\exists x(R \wedge T), \neg R, \exists
          x \neg T}}{ \root{\sruleRE}{\sqn{R \wedge T, \neg R, \exists
            x \neg T}}{ \rroot{\sruleRand}{\sqn{R \wedge T, \neg R,
              \neg T}}{ \root{\sruleAx}{\sqn{R,\neg R}}{ \leaf{}}}{
            \root{\sruleAx}{\sqn{T, \neg T}}{ \leaf{}}}}}}} } \quad ,
$$  
and for the rules 
$$
       \vcinf{\rulendown}
          {S\cons{\exists x R}}
          {S\cons{R\subst{x/t}}}
\qquad \mbox{and} \qquad
        \vcinf{\rulenup}
          {S\cons{R\subst{x/t}}}
          {S\cons{\forall x R}}
$$
we have the proofs
$$
\vcenter{\dernote{\sruleRE}{}{\sqn{\exists x R,
    \overline{R\subst{x/t}}}}{\root{\sruleAx}{
    \sqn{R\subst{x/t},\overline{R\subst{x/t}}} }{\leaf{}}}}
 \qquad\mbox{and} \qquad 
\vcenter{\dernote{\sruleRE}{}{\sqn{R\subst{x/t}, \exists x
    \neg R}}{\root{\sruleAx}{ \sqn{R\subst{x/t}, \neg R\subst{x/t}}
  }{\leaf{}}}}
\qquad .
$$
\end{proof}

\begin{Corollary}\label{cor:skstogs1}
  If a structure $S$ has a proof in $\SKSgpred$ then $\sqn{\lg{S}}$
  has a proof in $\GSone$.
\end{Corollary}

Soundness and completeness of $\SKSgpred$, \ie the fact that a
structure has a proof in $\SKSgpred$ if and only if it is valid,
follows from soundness and completeness of $\GSone$ by Corollaries
\ref{cor:gs1tosks} and \ref{cor:skstogs1}. Moreover, a structure $T$
implies a structure $R$ if and only if there is a derivation from $T$
to $R$, which follows from soundness and completeness and the
following theorem.

\begin{Theorem}\label{thm:dedpred}\theoremnl
  There is a derivation $\vcenter{\xy\xygraph{[]!{0;<1pc,0pc>:} {T}
      -@{=}^<>(.5){\strut\SKSgpred}[ddd] {R} }\endxy}$ if and only if
  there is a proof $\vcenter{\xy\xygraph{[]!{0;<1pc,0pc>:}
      {}*=<0pt>{}:@{|=}^<>(.5){\strut\SKSgpred}[ddd] {\pars{\neg T,R}}
    }\endxy}\quad .$
\end{Theorem}
\begin{proof}
  (same as the proof of Theorem \ref{thm:ded}.) 
  A proof $\Pi$ can be obtained from a given derivation $\Delta$ and a
  derivation $\Delta$ from a given proof $\Pi$, respectively, as follows:
  $$\vcenter{\xy\xygraph{[]!{0;<1pc,0pc>:}
      {\inf{\ruleidown}{\pars{\neg T, T}}{\true}\hspace{2.7ex}}
      -@{=}^<>(.5){\strut\SKSgpred} _<>(.5){\strut\pars{\neg
          T,\Delta}}[dddd] {\pars{ \neg T, R}} }\endxy} \qquad
  \mbox{and} \qquad \vcenter{\xy\xygraph{[]!{0;<1pc,0pc>:} {T}
      -@{=}_<>(.5){\aprs{T,\Pi}}^<>(.5){\strut\SKSgpred} [ddddd]
      {\dernote{\aintr}{}{R}{ \root{\swir}{\pars{R,\aprs{T,\neg T}}}{
            \leaf{\aprs{T,\pars{\neg T,R}} }}}} }\endxy} \qquad .
$$

\end{proof}

Note that the above is not true for system $\GSone$, because the
premise of the $\sruleRA$ rule does not imply its conclusion. The proof
above does not work for the sequent calculus because adding to the
context of a derivation can violate the proviso of the $\sruleRA$
rule.

\subsection{Admissibility of the Cut and the Other Up-Rules}

Just like in the propositional case, the up-rules of the symmetric
system are admissible.  By removing them from $\SKSgpred$ we obtain
the asymmetric, cut-free system shown in Figure~\ref{fig:ksgq}, which
is called system $\KSgpred$.

\renewcommand{\rulebox}[1]{\mbox{$#1$}}
\begin{figure}[htb]
  \begin{center}
    \fbox{
      \parbox{.8\textwidth}{
\medskip
$$
        \rulebox{\vcinf{\ruleidown}
          {S\pars{R,\neg R}}
          {S\cons{\true}}}
        \qquad\qquad
        \rulebox{\vcinf{\rulewdown}
          {S\cons{R}}
          {S\cons{\false}}}
        \qquad\qquad
        \rulebox{\vcinf{\rulecdown}
          {S\cons{ R}}
          {S\pars{ R,R }}}
$$
\medskip
$$
        \rulebox{\vcinf{\swir}
        {S\pars{\aprs{R,U},T }}
        {S\aprs{\pars{R,T},U}}}
\qquad
      \rulebox{\vcinf{\ruleudown}
          {S\pars{\forall x R, \exists x T}}
          {S\cons{\forall x \pars{R,T}}}}
\qquad
        \rulebox{\vcinf{\rulendown}
          {S\cons{\exists x R}}
          {S\cons{R\subst{x/t}}}}
$$
        }
      }    
    \caption{System $\KSgpred$}
    \label{fig:ksgq}
  \end{center}
\end{figure}

\renewcommand{\SKSproof}{\xy\xygraph{[]!{0;<1pc,0pc>:}
    {}*=<0pt>{}:@{|=}^<>(.5){\strut\SKSgpred}[dd] {S} }\endxy}
\renewcommand{\KSproof}{\xy\xygraph{[]!{0;<1pc,0pc>:}
    {}*=<0pt>{}:@{|=}^<>(.5){\strut\KSgpred}[dd] {S} }\endxy}
\renewcommand{\GScutproof}{
  \Derivation{\quad}{\hspace{9ex}}{\quad}{\sqn{\lg{S}}}{{\GSone \atop
+\sruleCut}}}
\renewcommand{\GSproof}{
  \Derivation{\quad}{\hspace{9ex}}{\quad}{\sqn{\lg{S}}}{\GSone}}

\begin{Theorem}\label{thm:transcutelimpred}
  The rules $\ruleiup$, $\rulewup$, $\rulecup$, $\ruleuup$ and
  $\rulenup$ are admissible for system $\KSgpred$.
\end{Theorem}

\begin{proof}
\theoremnl
\xymatrix@C11ex{
{\SKSproof} \ar[r]^-{\mbox{\scriptsize Corollary \ref{thm:SKSpredtoGS1}}} &  
{\GScutproof} \ar[r]^-{{\mbox{\scriptsize Cut
elimination} \atop \mbox{\scriptsize for $\GSone$}}} & 
{\GSproof} \ar[r]^-{\mbox{\scriptsize Corollary \ref{cor:gs1tosks2} }} & 
{\KSproof} \\
}
\end{proof}

\begin{Corollary}
  The systems $\SKSgpred$ and $\KSgpred$ are equivalent.
\end{Corollary}

\subsection{Reducing Rules to Atomic Form}

Consider the following local rules:
$$
\renewcommand{\arraycolsep}{2ex}
\begin{array}[h]{ll}
\inf{\rulemonedown} {S\cons{\exists
          x \pars{R,T}}} {S\pars{\exists x R, \exists x T}}
&
\inf{\rulemoneup} {S\aprs{\forall x R, \forall x T}}
      {S\cons{\forall x \aprs{R,T}}}
\\[2ex]
\inf{\rulemtwodown} {S\cons{\forall x \pars{R,T}}}
      {S\pars{\forall x R, \forall x T}}
&
\infnote{\rulemtwoup} {S\aprs{\exists x R, \exists x T}}
      {S\cons{\exists x \aprs{R,T}}}{\qquad .}
\end{array}
$$

Like medial, they have no analogues in the sequent calculus. In system
$\SKSgpred$, and similarly in the sequent calculus, the corresponding
inferences are made using contraction and weakening:

\begin{Proposition}\label{prop:predmedial}
  The rules $\{\rulemonedown,\rulemtwodown\}$ are derivable for
  $\{\rulecdown,\rulewdown\}$.  Dually, the rules $\{\rulemoneup,
  \rulemtwoup\}$ are derivable for $\{\rulecup,\rulewup\}$.
\end{Proposition}

\begin{proof}
  We show the case for $\rulemonedown$, the other cases are similar or
  dual:
  $$
  \dernote {\rulecdown }{\quad.}{S\cons{\exists x \pars{R,T}}} {
    \root{\rulewdown}{ S\pars{\exists x \pars{R,T}, \exists x
        \pars{R,T}}}{ \root{\rulewdown}{S\pars{\exists x R, \exists x
          \pars{R,T}} }{ \leaf{S\pars{\exists x R, \exists x T} } }}}
  $$
\end{proof}

Using medial and the rules $\{\rulemonedown,\rulemtwodown,\rulemoneup,
\rulemtwoup\}$ we can reduce identity, cut and weakening to atomic
form, similarly to the propositional case.

\begin{Theorem}\label{thm:predgeneral}
  The rules $\ruleidown$, $\rulewdown$ and $\rulecdown$ are derivable
  for $\{\ruleaidown,\swir, \ruleudown\}$, $\{\ruleawdown\}$ and
  $\{\ruleacdown, \rulem, \rulemonedown, \rulemtwodown \}$,
  respectively. Dually, the rules $\ruleiup$, $\rulewup$ and
  $\rulecup$ are derivable for $\{\ruleaiup,\swir, \ruleuup\}$,
  $\{\ruleawup\}$ and $\{\ruleacup, \rulem, \rulemoneup, \rulemtwoup
  \}$, respectively.
\end{Theorem}

\begin{proof}
  
  The proof is an extension of the proof of Theorem \ref{thm:general}
  by the inductive cases for the quantifiers.  Given an instance of
  one of the following rules:
$$
  \vcinf{\intr}{S{\pars{R,\neg R}}}{S\cons{\true}}\quad,\qquad
  \vcinf{\rulewdown}{S{\cons{R}}}{S\cons{\false}}\quad,\qquad
  \vcinf{\rulecdown}{S{\cons{R}}}{S\pars{R,R}} \quad,\qquad
  $$
construct a new derivation by structural induction on $R$:
\begin{enumerate}
\item $R=\exists x T$, where $x$ occurs free in $T$. Apply the
  induction hypothesis respectively on

$$\downsmash{
\dernote   {\ruleudown}{\quad,}  {S\pars{ \exists x T, \forall x \neg T }} {
\root   {\intr}  {S\cons{\forall x \pars{T,\neg T}} }{
\root   {=}  {S\cons{\forall x \tinyspace \true} }{
\leaf            {S\cons{\true}}                                   }}}
}\qquad
\downsmash{
\dernote   {\rulewdown}{\quad,}  {S\cons{\exists x T} }{
\rootnote   {=}{}  {S \cons{\exists x \tinyspace \false} }{
\leaf            {S\cons{\false}}}}                                   
}\qquad
\downsmash{
\dernote {\rulecdown}{\quad.}{S\cons{\exists x T}}  {
\root{\rulemonedown}{S\cons{\exists x \pars{T,T}}} {
\leaf                             {S\pars{\exists x T, \exists x T}}}}
}
$$
\item $R=\forall x T$, where $x$ occurs free in $T$. Apply the
induction hypothesis respectively on 
$$\downsmash{ \dernote {\ruleudown}{\quad,} {S\pars{\forall x T,
      \exists x \neg T }} { \root {\intr} {S \cons{\forall x
        \pars{T,\neg T}} }{ \root {=} {S\cons{\forall x \tinyspace
          \true} }{ \leaf {S\cons{\true}}}}} }\qquad 
\downsmash{
  \dernote {\rulewdown}{\quad,} {S\cons{\forall x T} }{ \root {=}
{S\cons{\forall x \tinyspace \false} }{ \leaf {S\cons{\false}}}} }\qquad
\downsmash{ 
\dernote{\rulecdown }{\quad.}{S\cons{\forall x T} } {
    \root{\rulemtwodown}{S\cons{\forall x \pars{T,T}}}{
      \leaf{S\pars{\forall x T,\forall x T }}}} }$$
\end{enumerate}

\end{proof}

We now obtain system $\SKSpred$ from $\SKSgpred$ by restricting
identity, cut, weakening and contraction to atomic form and adding the
rules
$\{\rulem,\rulemonedown,\rulemtwodown,\rulemoneup,\rulemtwoup\}$. It
is shown in Figure~\ref{fig:sksq}.

\renewcommand{\rulebox}[1]{\parbox{8em}{$$#1$$}}
\begin{figure}[tb!]
  \begin{center}
    \fbox{
      \parbox{12cm}{
        \hfill
        \rulebox{\vcinf{\ruleaidown}
          {S\pars{a,\neg a}}
          {S\cons{\true}}}
        \qquad
        \rulebox{\vcinf{\ruleaiup}
          {S\cons{\false}}
          {S\aprs{a,\neg a}}}
        \hfill
        \hspace{0mm}
\\ %\vspace{0mm} \dotfill \vspace{0mm} \\
       \vspace{0mm}
       \hfill
       \rulebox{\vcinf{\swir}
        {S\pars{\aprs{R,T},U}}
        {S\aprs{\pars{R,U},T}}}
        \hfill
        \vspace{-10mm}
\\
       \rulebox{\vcinf{\ruleudown}
          {S\pars{\forall x R, \exists x T}}
          {S\cons{\forall x \pars{R,T}}}}
        \hfill
        \rulebox{\vcinf{\ruleuup}
          {S\cons{\exists x \aprs{R,T}}}
          {S\aprs{\exists x R, \forall x T}}}
\\
        \rulebox{\vcinf{\rulemonedown}
          {S\cons{\exists x \pars{R,T}}}
          {S\pars{\exists x R, \exists x T}}}
        \hfill
        \rulebox{\vcinf{\rulemoneup}
          {S\aprs{\forall x R, \forall x T}}
          {S\cons{\forall x \aprs{R,T}}}} \vspace{-10mm}
\\
        \vspace{0mm}
        \hfill
        \rulebox{\vcinf{\rulem}
        {S\aprs{\pars{R,T},\pars{U,V} }}
        {S\pars{\aprs{R,U},\aprs{T,V} }}}
         \hfill
         \vspace{-10mm}
\\
        \rulebox{\vcinf{\rulemtwodown}
          {S\cons{\forall x \pars{R,T}}}
          {S\pars{\forall x R, \forall x T}}}
        \hfill
        \rulebox{\vcinf{\rulemtwoup}
          {S\aprs{\exists x R, \exists x T}}
          {S\cons{\exists x \aprs{R,T}}}}  
\\
        \vspace{0mm}
        \hfill
        \rulebox{\vcinf{\ruleawdown}
          {S\cons{a}}
          {S\cons{\false}}}
        \qquad
        \rulebox{\vcinf{\ruleawup}
          {S\cons{\true}}
          {S\cons{a}}}
        \hfill 
        \vspace{0mm}
\\ 
        \vspace{0mm}
        \hfill
        \rulebox{\vcinf{\ruleacdown}
          {S\cons{ a}}
          {S\pars{ a,a }}}
        \qquad
        \rulebox{\vcinf{\ruleacup}
          {S\aprs{ a,a }}
          {S\cons{ a}}}
        \hfill  
        \vspace{0mm}
\\
        \rulebox{\vcinf{\rulendown}
          {S\cons{\exists x R}}
          {S\cons{R\subst{x/t}}}}
        \hfill
        \rulebox{\vcinf{\rulenup}
          {S\cons{R\subst{x/t}}}
          {S\cons{\forall x R}}}

         }
      }    
    \caption{System $\SKSpred$}
    \label{fig:sksq}
  \end{center}
\end{figure}

As in all the systems considered, the up-rules, \ie $\{\rulenup,
\ruleuup,\rulemoneup, \rulemtwoup\}$ are admissible. Hence, system
$\KSpred$, shown in Figure~\ref{fig:ksq}, is complete.

\begin{Theorem}\label{thm:equiv-pred}
  System $\SKSpred$ and system $\SKSgpred$ are strongly equivalent. Also,
  system $\KSpred$ and system $\KSgpred$ are strongly equivalent.
\end{Theorem}

\begin{proof}
  Derivations in $\SKSgpred$ are translated to derivations in
  $\SKSpred$ by Theorem \ref{thm:predgeneral}, and vice versa by
  Proposition \ref{prop:predmedial}. The same holds for $\KSgpred$ and
  $\KSpred$.
\end{proof}

Thus, all results obtained for system $\SKSgpred$ also hold for system
$\SKSpred$. As in the propositional case, I will freely use general
identity, cut, weakening and contraction to denote a corresponding
derivation in $\SKSpred$ according to Theorem \ref{thm:predgeneral}.

\renewcommand{\rulebox}[1]{\mbox{$#1$}}
\begin{figure}[tb]
  \begin{center}
    \fbox{
      \parbox{\textwidth}{
\medskip
$$
        \rulebox{\vcinf{\ruleaidown}
          {S\pars{a,\neg a}}
          {S\cons{\true}}}
        \qquad\qquad
        \rulebox{\vcinf{\ruleawdown}
          {S\cons{a}}
          {S\cons{\false}}}
        \qquad\qquad
        \rulebox{\vcinf{\ruleacdown}
          {S\cons{ a}}
          {S\pars{ a,a }}}
$$
\medskip
$$
        \rulebox{\vcinf{\swir}
        {S\pars{\aprs{R,U},T }}
        {S\aprs{\pars{R,T},U}}}
      \qquad
        \rulebox{\vcinf{\rulem}
        {S\aprs{\pars{R,U},\pars{T,V} }}
        {S\pars{\aprs{R,T},\aprs{U,V} }}}
$$
\medskip
$$
      \rulebox{\vcinf{\ruleudown}
          {S\pars{\forall x R, \exists x T}}
          {S\cons{\forall x \pars{R,T}}}}
\qquad
        \rulebox{\vcinf{\rulendown}
          {S\cons{\exists x R}}
          {S\cons{R\subst{x/t}}}}
\qquad
       \rulebox{\vcinf{\rulemonedown}
          {S\cons{\exists x \pars{R,T}}}
          {S\pars{\exists x R, \exists x T}}}
\qquad
        \rulebox{\vcinf{\rulemtwodown}
          {S\cons{\forall x \pars{R,T}}}
          {S\pars{\forall x R, \forall x T}}}
$$
        }
      }    
    \caption{System $\KSpred$}
    \label{fig:ksq}
  \end{center}
\end{figure}

\subsection{Locality Through Atomicity}

As we have seen in the previous section, the technique of reducing
contraction to atomic form to obtain locality also works in the case
of predicate logic: the non-local rule $\rulecdown$ is equivalently
replaced by local ones, namely
$\{\ruleacdown,\rulem,\rulemonedown,\rulemtwodown\}$.

However, there are other sources of non-locality in system
$\SKSpred$. One is the condition on the quantifier equations:
$$
\begin{array}{l}
\forall y R = \exists y R = R \qquad \mbox{where $y$ is
not free in $R$.}
\end{array}
$$
To add or remove a quantifier, a structure of unbounded size has to
be checked for occurrences of the variable $y$.  

Another is the $\rulendown$ rule, in which a term $t$ of unbounded
size has to be copied into an unbounded number of occurrences of $x$
in $R$. It is non-local for two distinct reasons: 1) the unbounded
size of $t$ and 2) the unbounded number of occurrences of $x$ in
$R$. The unboundedness of term $t$ can be dealt with, since
$\rulendown$ can be derived and thus replaced by the following two
rules:
$$
\infnote{{\sf n_1}{\downarrow}}{S\cons{\exists x R}}{S\cons{\exists
y_1 \dots \exists y_n R\subst{x/f(y_1,...,y_n)}}}{\qquad \mbox{and}
\quad } \qquad \infnote{{\sf
    n_2}{\downarrow}}{S\cons{\exists x R}}{S\cons{R}}{\quad ,}
$$
where $f$ is a function symbol of arity $n$. Still, rule ${\sf
  n_1}{\downarrow}$ is not local because of the unbounded number of
occurrences of $x$ in $R$.

Is it possible to obtain a local system for first-order predicate
logic?  I do not know how to do it without adding new symbols to the
language of predicate logic. But it is conceivable to obtain a local
system by introducing substitution operators together with rules that
explicitly handle the instantiation of variables piece by piece. The
question is whether this can be done without losing the good
properties, especially cut elimination and simplicity.

\subsection{Decomposition of Derivations}

In the following I show how all decomposition results for the
propositional system from Section \ref{sec:decomp} extend to predicate
logic in a straightforward way.

As in the propositional case, atomic identity and cut can be reduced
to their shallow versions using the super switch rules. In the
predicate case the rules shallow atomic identity and shallow atomic
are as follows:
$$
\vcinf{\ruleaisdown} {\aprs{S,\forall \pars{a,\neg a}}} {S} \qquad{\rm and}\qquad \vcinf{\ruleaisup} {S}
{\pars{S,\exists \aprs{a,\neg a}}} \quad ,
$$
where $\forall$ and $\exists$ denote sequences of quantifiers that
universally close $\pars{a,\neg a}$ and existentially close
$\aprs{a,\neg a}$, respectively.

The super switch rules for predicate logic, 
$$
\vcinf{\rulessdown} {S\pars{R,T\cons{\false}}} {S\cons{T\cons{R}}}
\qquad{\rm and}\qquad \vcinf{\rulessup} {S\cons{T\cons{R}}}
{S\aprs{R,T\cons{\true}}} \quad ,
$$
carry a proviso: quantifiers in $T$ do not
 capture variables in
$R$. This is not a restriction because bound variables can always be
renamed such that the proviso is fulfilled.

\begin{Lemma}\label{lem:superswitch-pred}
  The rule $\rulessdown$ is derivable for $\{\swir, \rulendown,
  \ruleudown\}$. Dually, the rule $\rulessup$ is derivable for $\{\swir,
  \rulenup, \ruleuup\}$.
\end{Lemma}
\begin{proof}
  The proof is an extension of the proof of Lemma
  \ref{lem:superswitch}.    I show the
  two cases that have to be considered in addition to the proof in the
  propositional case:
\begin{enumerate}
\item $T\cons{\enspace}=\forall x U\cons{\enspace}$, where
  $x$ occurs freely in $U$. Apply the induction hypothesis on
$$\downsmash{ 
\dernote {\rulessup }{}{S\cons{\forall x U\cons{R}}} {
\root{\rulenup }{S\cons{\forall x \aprs{R,U\cons{\true}}}} {
     \root{=}{S\cons{\forall x \aprs{R,\forall x U\cons{\true}}}}{
         \leaf{S\aprs{R,\forall x U\cons{\true}} }}}} }
\quad .$$
\item  $T\cons{\enspace}=\exists x U\cons{\enspace}$, where
  $x$ occurs freely in $U$. Apply the induction hypothesis on
  $$\downsmash{ \dernote {\rulessup }{}{S\cons{\exists x U\cons{R}}} {
      \root{\ruleuup }{S\cons{\exists x \aprs{R, U\cons{\true}}}} {
        \root{=}{S\aprs{\forall x R,\exists x U\cons{\true}}}{
          \leaf{S\aprs{R,\exists x U\cons{\true}} }}}} } \quad .$$
\end{enumerate}
\end{proof}

\begin{Lemma}\label{lem:shallowint-pred}
  The rule $\ruleaidown$ is derivable for $\{\ruleaisdown,\swir,
  \rulenup, \ruleuup\}$. Dually, the rule $\ruleaiup$ is derivable for
  $\{\ruleaisup,\swir, \rulendown, \ruleudown\}$.
\end{Lemma}
\begin{proof}
 
$$
\mbox{An instance of} \quad
\vcinf{\ruleaidown}{S\pars{a,\neg a}}{S\cons{\true}}
\quad \mbox{is replaced by} \quad 
 \vcenter{\dernote{\rulenup^n}{}{S\pars{a,\neg a}}{
    \root{\rulessup}{S\cons{\forall \pars{a,\neg a}}}{
    \root{\ruleaidown}{\aprs{S\cons{\true},\forall \pars{a,\neg a}}}{
      \leaf{S\cons{\true}}}}}}\quad .
$$
(And dually for $\ruleaiup$.)
\end{proof}

\begin{Theorem}[Decomposition in Predicate Logic]\label{thm:decomp-pred}
  All theorems of section \ref{sec:decomp} also hold in the case of
  predicate logic, \ie with $\SKS$ replaced by $\SKSpred$ and $\KS$
  replaced by $\KSpred$.  In Theorem \ref{thm:sepall},
  $\{\rules,\rulem\}$ has to be extended by the quantifier rules, \ie
  $\{\ruleudown,\ruleuup,\rulemonedown,\rulemoneup,\rulemtwodown,\rulemtwoup,\rulendown,\rulenup\}$.
\end{Theorem}

\begin{proof}
  Identity and cut are separated as in the propositional case, using
  Lemma \ref{lem:shallowint-pred} instead of Lemma
  \ref{lem:shallowint}.
  
  Contraction is separated as in the propositional case, using the
  proof of Theorem \ref{thm:sepac2}. The only difference is in step
  four, where instances of $\ruleacdown$ have to be permuted under
  instances of rules from $\KSpred\setminus\KS$. None of those rules
  except for $\rulendown$ changes atoms, so $\ruleacdown$ trivially
  permutes under those instances. It also easily permutes under
  instances of $\rulendown$:
$$
  \dernote{\rulendown}{\quad \leadsto
    \quad}{S\cons{\exists x R\cons{a}}}{
    \root{\ruleacdown}{S\cons{R\cons{a}\subst{x/t}}}{
      \leaf{S\cons{R\pars{a,a}\subst{x/t}}}}}
  \dernote{\ruleacdown}{\quad .}{S\cons{\exists x R\cons{a}}}{
    \root{\rulendown}{S\cons{\exists x R\pars{a,a}}}{
      \leaf{S\cons{R\pars{a,a}\subst{x/t}}}}}
$$

Weakening is separated as in the propositional case. When moved over
$\rulendown$ and $\rulenup$, derivations of weakenings will contain
weakenings on different atoms (with the same predicate symbol but
differently instantiated):
$$
  \dernote{\ruleawup^n}{\quad \leadsto
    \quad}{S\cons{\exists x R\cons{\true}}}{
    \root{\rulendown}{S'\cons{\exists x R'\cons{a}}}{
      \leaf{S'\cons{R'\cons{a}\subst{x/t}}}}}
  \dernote{\rulendown}{}{S\cons{\exists x R\cons{\true}}}{
    \root{\ruleawup^n}{S\cons{R\cons{t}\subst{x/t}}}{
      \leaf{S'\cons{R'\cons{a}\subst{x/t}}}}}
$$
\end{proof}

\section{Conclusions}

We have seen deductive systems for classical propositional and
predicate logic in the calculus of structures. They are sound and
complete, and the cut rule is admissible. In contrast to sequent
systems, their rules apply \emph{deep} inside formulae, and
derivations enjoy a top-down \emph{symmetry} which allows to dualise
them.

Those features allow to reduce the cut, weakening and contraction to
atomic form, which is not possible in the sequent calculus. This leads
to \emph{local} rules, \ie rules that do not require the inspection of
expressions of unbounded size. For propositional logic, I presented
system $\SKS$, which is local, \ie contains only local rules. For
predicate logic I presented system $\SKSpred$ which is local except for
the treatment of variables.

The freedom in applying inference rules in the calculus of structures
allows permutations that can not be observed in the sequent
calculus.  This leads to more normal forms for derivations, as shown
in the decomposition theorems.

Normal forms for derivations are an interesting area for future work.
The decomposition theorems given here barely scratch the surface of
what seems to be achievable.  For example, consider the following two
conjectures:

\begin{Conjecture}[Interpolation] \theoremnl
   For every derivation $\vcenter{\xy\xygraph{[]!{0;<1pc,0pc>:} {T}
        -@{=}^<>(.5){\strut\SKS}[dd] {R}
        }\endxy}$  there is a derivation
    $\vcenter{\xy\xygraph{[]!{0;<1pc,0pc>:} {T}
      -@{=}^<>(.5){\strut\SKS\setminus\{\ruleaidown,\ruleawdown\}}[dd] {P}    
    -@{=}^<>(.5){\strut\SKS\setminus\{\ruleaiup,\ruleawup\}} [dd] {R}
}\endxy}$.

\end{Conjecture}

Here, a derivation is separated into two phases: the top one, with
rules that do not introduce new atoms going down, and the bottom one,
with rules that do not introduce new atoms going up. Consequently, the
structure $P$ contains only atoms that occur in both $T$ and $R$ and
is thus an \emph{interpolant}. This form of interpolation can be seen
as the symmetric closure of cut elimination: not only are cuts pushed
up, but also their duals, identities, are pushed down. Cut elimination
is an immediate corollary of this property: if $T$ is equivalent to
the unit true then also $P$ is equivalent to true, and in the bottom
part of the derivation there are no cuts.  I will show elsewhere a
semantic proof of interpolation in the propositional case. A syntactic
proof that scales to the predicate case would be desirable.

Another decomposition theorem, that has been proved for two other
systems \cite{GugStr01} in the calculus of structures and led to cut
elimination, is the separation of the \emph{core} and the
\emph{non-core} fragment. So far, all the systems in the calculus of
structures allow for an easy reduction of both cut and identity to
atomic form by means of rules that can be obtained in a uniform way.
Those rules are called the \emph{core} fragment. In $\SKS$, the core
consists of one single rule: the switch. The core of $\SKSpred$, in
addition to the switch rule, also contains the rules $\ruleudown$ and 
$\ruleuup$.  All rules that are not in the core and are not identity
or cut are called \emph{non-core}. The problem is separating switch
and medial.

\begin{Conjecture}[Separation core -- non-core]\label{con:fullsep}
\theoremnl
    For every derivation $\vcenter{\xy\xygraph{[]!{0;<1pc,0pc>:} {T}
        -@{=}^<>(.5){\strut\SKS}[dd] {R}
        }\endxy}$  there is a derivation
    $\vcenter{\xy\xygraph{[]!{0;<1pc,0pc>:} {T}
    -@{=}^<>(.5){\strut\mbox{\emph{non-core}}}[dd] {T'}
    -@{=}^<>(.5){\strut\{\ruleaidown\}}[dd] {T''}
    -@{=}^<>(.5){\strut\mbox{\emph{core}} }[dd] {R''}
    -@{=}^<>(.5){\strut\{\ruleaiup\}}[dd] {R'}    
    -@{=}^<>(.5){\strut\mbox{\emph{non-core}}} [dd] {R}
}\endxy}$.
    
\end{Conjecture}

A cut elimination procedure that is based on permuting up instances of
the cut would be easy to obtain, could we rely on this conjecture.
Then all the problematic rules that could stand in the way of the cut
can be moved either below all the cuts or to the top of the proof,
rendering them trivial, since their premise is \emph{true}.  Cut
elimination is thus an easy consequence of such a decomposition
theorem.

The proof of the separation of contraction (Theorem \ref{thm:sepac2})
relies on the admissibility of the cut. It should be provable
directly, \ie without using cut admissibility, just by very natural
permutations. The difficulty is in proving termination of the process
of bouncing contractions up and down between cuts and identities, as
happens in \cite{StraELS01}.

The above mentioned freedom in applying inference rules is a mixed
blessing. Compared to the sequent calculus, it implies a greater
non-determinism in proof search. It will be interesting to see whether
it is possible to restrict this non-determinism by finding a suitable
notion of \emph{uniform proofs} \cite{MilNadPfeSceAPAL91}. Another
interesting question for future research is whether there is a local
system for intuitionistic logic.

\noindent
{\bf Acknowledgements}\\
This work has been accomplished while I was supported by the DFG
Gra\-du\-ier\-ten\-kol\-leg 334.  I would like to thank the members of the proof
theory group in Dresden for providing an inspiring environment,
especially Alessio Guglielmi. He also helped me with this paper in
numerous ways. Steffen Hölldobler, Lutz Straßburger and Charles Stewart
carefully read preliminary versions of this paper and made helpful
suggestions.

\noindent
{\bf Web Site}\\
Information about the calculus of structures is available from the
following URL: 
\begin{center}
  http://www.wv.inf.tu-dresden.de/\char'176 guglielm/Research/ .
\end{center}

\bibliographystyle{plain}
\bibliography{kai}

\end{document}